\pdfoutput=1
\documentclass[12pt,notitlepage]{article} 
\usepackage[heightrounded, marginparwidth=0.8in, marginparsep=0.15in]{geometry}  
\usepackage[small,nohug,heads=vee]{diagrams}
\usepackage{graphicx}
\usepackage{amssymb}
\usepackage{amsmath}
\usepackage{subcaption}
\usepackage[backend=bibtex]{biblatex}
\usepackage{float}
\bibliography{Automatic_Logarithm_for_arXiv_2018_11_30}

\newtheorem{theorem}{Theorem}[section]
\newtheorem{lemma}[theorem]{Lemma}
\newtheorem{remark}[theorem]{Remark}

\newtheorem{corollary}[theorem]{Corollary}
\newtheorem{proposition}[theorem]{Proposition}
\newtheorem{definition}[theorem]{Definition}
\newtheorem{example}[theorem]{Example}

\def \defeq {:=}
\def \sqr {\ \ensuremath{\square}}
\def \tri {\ \ensuremath{\triangle}}
\def \bb {\textbf}
\def \ii {\textit}
\def \pp {\par \noindent}
\def \pf {\pp \bb{Proof.\ }} 
 
\def\lp{\left(}
\def\rp{\right)}

\newcommand{\nn}{\mathbb{N}}
\newcommand{\zz}{\mathbb{Z}}
\newcommand{\rr}{\mathbb{R}}
\newcommand{\pr}{\mathbb{P}}

\def\<{\left\langle }
\def\>{\right\rangle}

\def\<-{\leftarrow }
\def\->{\rightarrow} 
\def\<={\Leftarrow }
\def\=>{\Rightarrow} 
\def\op{\operatorname}

\newcommand{\cyl}[1]{#1 X^\nn}	
\newcommand{\pcyl}[1]{\lp \cyl{#1} \rp} 
\newcommand{\alesh}{automaton $F$}
\newcommand{\zoran}{automaton $Z$}

\def\A{\mathcal{A}}

\def\O{\mathcal{O}}

\def\L{\mathcal{L}}

\def\mathcalA{\mathcal{A}}
\def\mathcalB{\mathcal{B}}
\def\mathcalC{\mathcal{C}}
\def\mathcalO{\mathcal{O}}
\def\mathcalS{\mathcal{S}}
\def\mathcalG{\mathcal{G}}
\def\mathcalT{\mathcal{T}}
\def\mathcalL{\mathcal{L}}

\def\mathbbC{\mathbb{C}}

\def\sem{\operatorname{sem}}
\def\gr{\operatorname{gr}}

\def\Log{\operatorname{Log}}
\def\Exp{\operatorname{Exp}}
\title{Automatic Logarithm and Associated Measures}
\author{Rostislav Grigorchuk \thanks{The author was supported by NSA grant H98230-15-1-0328 
and by Simons Foundation Collaboration Grant for Mathematicians, Award Number 527814.} \and Roman Kogan \and Yaroslav Vorobets}

\begin{document}

\maketitle
\begin{abstract}
We introduce the notion of the \emph{automatic logarithm} $\mathcal L_{\mathcal A, \mathcal B}$ with the purpose of studying the expanding properties of Schreier graphs of action of the group generated by two finite initial Mealy automata $\mathcal A$ and  $\mathcal B$ on the levels of a regular $d$-ary rooted tree $\mathcal T$, where $\mathcal A$ is level-transitive and of bounded activity. $\mathcal L_{\mathcal A, \mathcal B}$ computes the lengths of chords in this family of graphs. Formally, $\mathcal L$ is a map $\partial \mathcal T \rightarrow \zz_d$ from the boundary of the tree to the integer $p$-adics whose values are determined by a Moore machine. The distribution of its outputs yields a probabilistic measure $\mu$ on $\partial \mathcal T$, which in some cases can be computed by a Mealy-type machine  (we then say that $\mu$ is \emph{finite-state}). We provide a criterion to determine whether $\mu$ is finite-state. A number of examples illustrating the different cases with $\mathcal A$ being the adding machine is provided.
\end{abstract}

\section{Introduction}

The maps and the measures considered in this paper arise from the study of properties of Schreier graphs  associated with automaton semigroups and groups acting on words over a finite alphabet and regular rooted  trees.

The problem of studying the distribution of lengths of chords (to be defined below) in the graph of action of two initial automata gives rise to the \ii{automatic logarithm}, a map defined by an automaton that outputs these lengths. The distribution of the lengths of chords is then seen as the image of the uniform Bernoulli measure by the action of the automatic logarithm. When the automatic logarithm is invertible, the distribution is uniform. Otherwise, the resulting distribution is an interesting object of study. In certain cases, such distributions only have a finite number of restrictions to cylinders (we call them \bb{finite-state}), and we provide a sufficient condition for this to happen, as well as examples when it does not.

Given  a  finite initial Mealy (or Moore) type  automaton $\mathcalA_q$ over a finite  alphabet $X$, one can define a  map  $\hat{\mathcalA_q}$ on the  space of sequences (words) over the alphabet  $X$. Maps of this type usually  have a very complicated dynamical nature and may  transform relatively  simple  measures  on the  space $X^{\mathbb{N}}$,  like for instance  Bernoulli  or   Markov measure, into complicated ones.  The  study  of  such  measures were  initiated in  \cite{AKP}, \cite{Ryab} and \cite{Krav}.

Given  a family  of  finite  automata  $\mathcalA_q, \mathcalB_s,\dots, \mathcalC_t$,  using  the  operation  of  composition  of  automata one can generate  a semigroup $\mathcalS= \langle\mathcalA_q, \mathcalB_s,\dots, \mathcalC_t\rangle_{\sem}$  or  even a  group  $\mathcalG=\langle\mathcalA_q, \mathcalB_s,\dots, \mathcalC_t\rangle_{\gr}$  if  the  automata are  invertible.  A  particularly  interesting  case  is   when  the  group $\mathcalG$  is  generated  by  a  family  which  comes  from  the one  non-initial invertible  automaton  $\mathcalA$  by  using  all  its  states for  generating.  Such  groups  are called automaton groups  (or  self-similar  groups)  and  play  an important  role  in group  theory  and  areas  of  its  applications  \cite{GB}, \cite{Nek}, \cite{GNS}, \cite{GNSunic}.   They  naturally  act  by  automorphisms  on a $d$-regular  rooted  tree $\mathcalT$  ($d$  is  the  cardinality  of  the  alphabet  $X$)  and  on  its  boundary $\partial T$.  These  actions  are   induced  by  the  corresponding  actions  on  the  set  of  finite (and, respectively, infinite)  words over the  alphabet  $X$.  The  operation  of  composition  of  automata  corresponds  to  the  composition  of the associated  maps.

Another direction  of  development  is  study  of the  Schreier  graphs (also called orbital  graphs)  given  by  the  action  of  a  group  on  levels of the  tree   or  on  its  boundary (i.e.   on  finite  or  infinite  words). These  graphs have  self-similarity  features and give a  good  approximation  to many  important  fractal  sets   including the Julia  sets of  the  rational  mappings  of    $\mathbbC$.   There  are   examples  of  the  automata given  by  a  small  number  of  states that are believed to produce families of  expander  graphs (two  of  them  are  considered  in this  article). No  rigorous  proof  of  this is  known, but  there are   results  showing  that  at  least  these families  are the so  called  asymptotic  expanders  \cite{GrigExp},  and  that  the  growth  of  their  diameters   is  slow \cite{PakMal} (as  should  be  in  the  case  of  expanders).

Among   automorphisms  of  the  rooted  trees,  the most famous is the \ii{adding machine} automorphism defined by the automaton shown in Figure \ref{fig:odom1}, which we denote $\mathcalO$. The \ii{portrait} of this automorphism is shown in Figure \ref{fig:odomportrait}. It acts on finite  strings  of  symbols as the addition of $1$ when the strings are interepreted as the natural  numbers in the $d$-adic expansion (the diagram in Figure \ref{fig:odom1} is for the case of the binary alphabet). The  group  $\mathcal{G}(\mathcalO)$ is an infinite cyclic group (one of the states of $\mathcalO$ corresponds to the identity map). If $o$ is the nontrivial state of $\mathcalO$, and if another initial automaton $\mathcalA_a$ is given, then one can consider the semigroup $\langle A_a,  O_o\rangle_{\sem}$ (or  a  group $\langle \A_a,  \O_o\rangle_{\gr}$ if  $\mathcalA_a$  is  invertible), and  study  its properties  as  well as  the  sequence  $\{\Gamma_n\}_{n=1}^{\infty}$ of graphs of action on the  levels  $n=1,2,\dots,n\dots$ of  the tree  
(Figure  \ref{fig:schreiergraphexamples}   gives an  impression  of  how  the  graph $\Gamma_n$ may  look). The  questions  about   combinatorial  and  spectral  properties  of  graphs $\{\Gamma_n\}$ is the subject of  many  investigations \cite{GBspec, GZ, GSHanoi}, and in particular, the  question about the growth of the diameters of $\{\Gamma_n\}$  and  about the expansion properties of the family $\{\Gamma_n\}_{n=1}^{\infty}$   are   among  the  central.

In  this  paper  we  focus  on  study of  the  dynamical and combinatorial properties of  the  pair ($\mathcalA_a, \mathcalO_o$). This  unexpectedly  leads  us  to  the  notion  of  the  ``logarithm''  of  $\mathcalA_a$  with  respect  to  the adding machine $\mathcalO_o$, which  we  denote $\mathcal L$, and which is a specially defined map  on  the  set  of  finite  and  infinite  sequences  with vales in $d$-adic numbers, or words over $X$ that represent them. Then we show that in the case of the binary alphabet, the logarithm  map can be define by a finite state  automaton  (Theorem  \ref{dautom}) and provide the construction for it.  We then analyze the distribution of the lengths  of  the  ``chords''  
(again we appeal to Figure \ref{fig:schreiergraphexamples} which  gives an impression of what   we  mean   by  the  chord).  This  leads  us  to  the   considerations started  in  the \cite{Krav, GKV1} about the  nature  of  the  image  of  the  Bernoulli  (or, more  generally, Markov)  measure  under  the  automaton map,  in  the  case  the  map  is  given  by  the  
``logarithm''  automaton  $\mathcalL$.  The  distribution  of the chords is given  by the image $\mu = \mathcal{L}_{*}(\nu)$ of the uniform Bernoulli  measure  $\nu$ on  $X^\nn$, which  in  some  important  cases  (for  instance,
given  in the Example \ref{bellaex} and Theorem \ref{aleshmex})  is a  Markov  measure, but in some other interesting cases  (like  the  Example \ref{lampex}) is a more  complicated type of measure.  

To study $\mu$, we  introduce  the  notion  of  the  automaton associated  with  a measure, the  notion of a   finite  state  measure, of a self-similar  measure,   and  show  that Markov  measure  is  a finite  state  measure,  but  not  every  finite  state  measure   is  Markov. Multiple examples illustrating the definitions and results are provided (we mark the ends of examples with $\tri$ instead of the traditional $\sqr$).

\section{Preliminaries}
\subsection{Endomorphisms of rooted trees}

By a finite \bb{alphabet} $X$ of size $d$, we mean  a finite set of cardinality $d$. We will usually use the following alphabets: the sets $\{0,1, \ldots, d-1\}$, and finite subsets $X \subset [0,1]$ (in the latter case, the values of the real numbers constituting $X$ will be important in addition to its cardinality).

For a word $w$ in $X$, $|w|$ denotes its length, and $w_i$ denotes the $i$'th character for $0\leq i \leq |w|-1$.
The numbering of characters starts from $0$, so $w=w_0w_1\ldots w_{|w|-1}$. If $v$ is another word (or a character), $wv$ is the concatenation of the two. 

$X^*$ denotes all finite words over $X$:
\[
X^* := \{ a_0\ldots a_{n-1}\ : \ a_i \in X, n \in \nn \cup \{0\}\}
\]

Let $\mathcal T$ be a rooted graph with the vertex set $V=X^*$, edge set $ F=\{(w, wa)\ :\ w \in X^*, a\in X\}$, and the root being the empty word. This graph is a \bb{$d$-regular rooted tree}.


The $n$'th level $X^n$ of the tree $\mathcal T$ is the set of words of length $n$.
 
An \bb{endomorphism} of the rooted tree $\mathcal T$ is a map $g$ from $X^*$ to itself that preserves the levels and maps adjacent vertices to adjacent vertices. An \bb{automorphism} is an invertible endomorphism.

The \bb{boundary} of the tree $\mathcal T$ is the set $X^\nn$ of infinite sequences in $X$:
\[
\partial \mathcal T := \{ a_0a_1a_2\ldots\ :\ a_i \in X, i \in \nn\}
\] 
$\partial T$ is supplied with the Tychonoff product topology that makes it homeomorphic to a Cantor set. Geometrically, the boundary can be viewed as a set of geodesic paths starting at the root and going to infinity.

Let $\sigma_r$ denote the operation that deletes the last character of a word: for $w \in X^*$ and $a \in X$,
$\sigma_r(wa) := w$.
Then $\partial T$ can be obtained as the inverse limit of the directed system of levels $\{X^n\ :\ n \in \nn\}$ with the projections $\psi_{m,n} : X^n \-> X^m$ given by $\psi_{m,n} := \sigma_r^{n-m}$ (i.e., discarding the last $n-m$ characters).

\subsection{Mealy and Moore machines}

\begin{figure}[ht]
\centering
	\begin{subfigure}{0.45\textwidth}
		\centering	
		\includegraphics[width=0.75\textwidth]{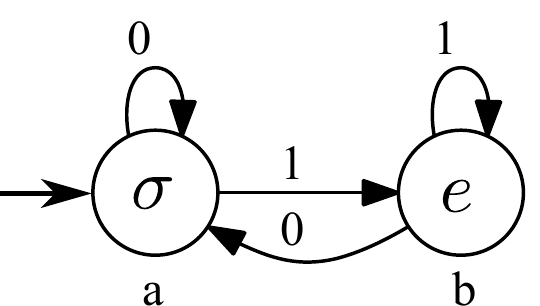}
		\caption{Diagram with element of the symmetric group on vertices}
		  		 \label{fig:lamp2}
	\end{subfigure}
		\begin{subfigure}{0.45\textwidth}
		\centering	
		\includegraphics[width=0.6\textwidth]{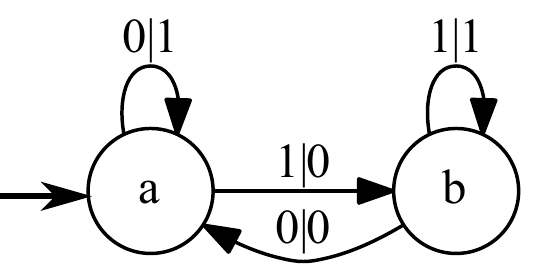}
		 \caption{Diagram with output marked on edges}
		 		\label{fig:lamp1}
	\end{subfigure}

	\caption{Two ways to draw the Lamplighter automaton}
	\label{fig:automdiagrams}
\end{figure}

\begin{definition}
A \bb{Mealy machine}, or a \bb{finite initial automaton with output}, is a hextuple $\A_q=(S, q, X, Y, \pi,\lambda)$, where
\begin{itemize}
\item $X$ is a (finite) input alphabet;
\item $Y$ is a (finite) output alphabet;
\item $S$ is a (finite) set of \ii{states};
\item $q \in S$ is the \ii{initial state};
\item $\pi: S \times X \-> S$ is the \ii{transition map}
\item $\lambda: S\times X \-> Y$ is the \ii{output map}
\end{itemize}

When the initial state of the automaton is understood from the context, we drop the subscript and write $\mathcal A$ instead of $\mathcal A_q$ to denote it.

We write $\pi_s$, $\lambda_s$ for restrictions of these functions to the state $s$, defining $\pi_s(x) := \pi(s,x)$ and $\lambda_s(x) :=\lambda(s,x)$.
 
The functions $\pi$ and $\lambda$ also act on words in the alphabet $X$ via the following recursive relations (for $x\in X$, $w \in X^*$):
\begin{align*}
\pi(s, xw) &:= \pi(\pi(s,x), w);\\
\lambda(s, xw) &:= \lambda(s,x)\lambda(\pi(s,x), w).\\
\end{align*}
In the same way, $\pi_s$ and $\lambda_s$, for $s \in S$, act on words $w\in X^*$. Additionally, we may write $\pi(w)$ for $\pi_q(w)$ (and similarly, $\lambda(w)$ for $\lambda_q(w)$), when the initial state $q$ is understood from the context.
\end{definition}
The \bb{diagram} of an automaton $\A_q$ is a labeled graph with the vertex set $S$, edge set $E=\{(s, \pi(s,x))\ :\ s \in S, x \in X\}$, with label $x:\lambda(s,x)$ on the edge $(s, \pi(s,x))$. 
The initial state $q$ is marked with a special arrow (which doesn't start at a state).
An example of such diagram for the Lamplighter automaton is shown in Figure \ref{fig:lamp1}.

An automaton $\mathcal A$ is \bb{invertible} if $\lambda_s$ is invertible for all $s \in S$ (that is, if $\lambda_s \in S(X)$, where $S(X)$ is the symmetric group on $X$). The endomorphism $g$ given by an invertible automaton $\mathcal A$ is invertible, and the automaton for $g^{-1}$ (which we denote as $\mathcal A^{-1}$) can be constructed from the diagram of $\mathcal A$ by flipping the input and output on the edges.

In the case when an automaton is invertible, we can draw the \bb{diagram} of the automaton without specifying its output on the arrows. Instead, the state $s$ is marked by the element of the symmetric group $\lambda_s \in S(X)$. If $\lambda_s$ is the trivial permutation, we call the state $s$ \bb{passive}, and call it \bb{active} otherwise.

\begin{figure}[ht]
\centering
	\begin{subfigure}{0.32\textwidth}
		\centering	
		\includegraphics[width=0.6\textwidth]{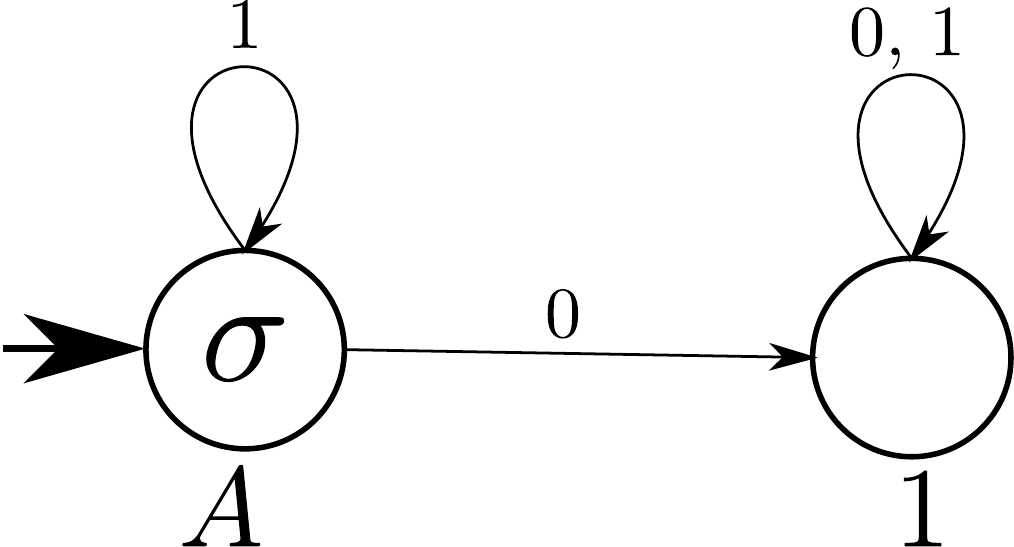}
		 \caption{Adding machine, also featured in Figure \ref{fig:odomportrait}}
		 		\label{fig:odom1}
	\end{subfigure}
	\begin{subfigure}{0.32\textwidth}
		\centering	
		\includegraphics[width=0.6\textwidth]{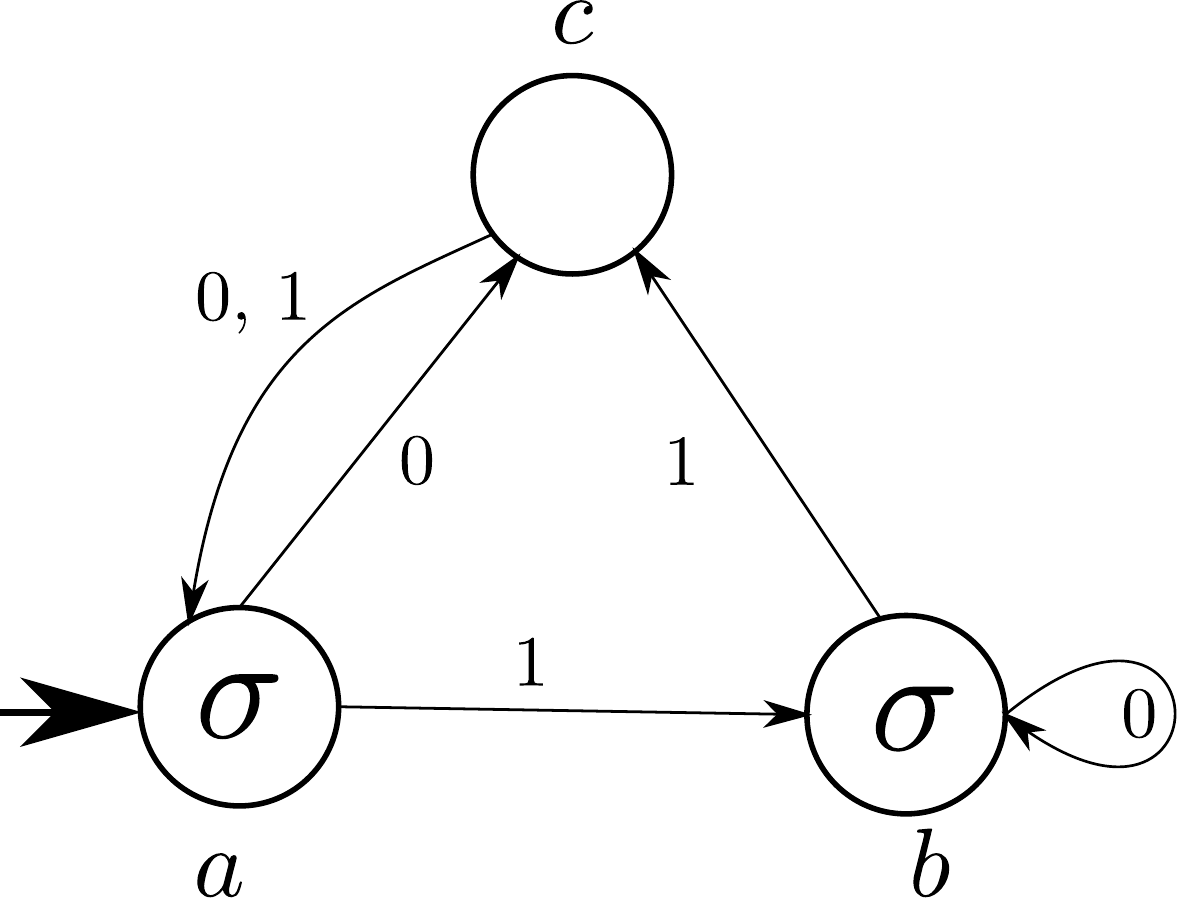}
		\caption{\alesh}
		  		 \label{fig:alesh1}
	\end{subfigure}
		\begin{subfigure}{0.32\textwidth}
		\centering	
		\includegraphics[width=0.6\textwidth]{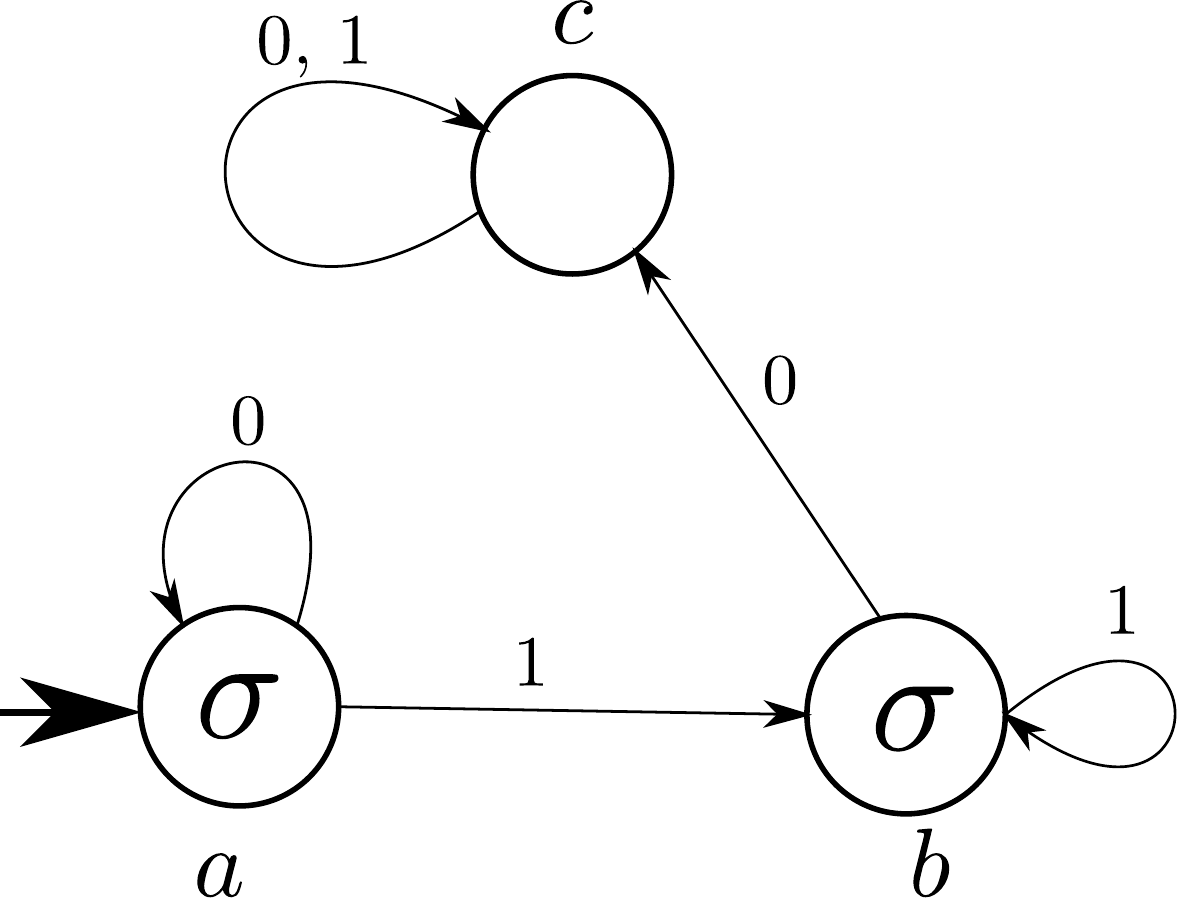}
		 \caption{\zoran }
		 		\label{fig:zoran}
	\end{subfigure}
	\caption{Diagrams of invertible automata}
	\label{fig:automexamples}
\end{figure}

When $X=\{0,1\}$, we write $\sigma$ for the nontrivial permutation of $X$ (i.e. $\sigma(0)=1, \sigma(1)=0$). In the diagrams of automata over $X=\{0,1\}$ we then mark active states with $\sigma$, leave the label of passive states blank. 
Figure \ref{fig:lamp2} shows how to draw the Lamplighter automaton of Figure \ref{fig:lamp1} in this way. A few more examples of such diagrams are in Figure \ref{fig:automexamples}, and further throughout this paper.

Unless otherwise specified, we assume $X=Y$ everywhere in this text, and write an automaton $\A_q=(S, q, X, \pi, \lambda)$. 

An automaton state $q$ acts on $X^*$, the $d$-ary tree $\mathcal T$, and its boundary $\partial T$ by the action of $\mathcal A_q$. We shall use $\mathcal A$ and $q$ interchangeably for this action when the context is clear. 

\begin{definition} The \bb{graph of the action} of an initial automaton $\A_q=(S,q,X,\pi,\lambda)$ on an invariant subset $S \subset V(\mathcal T)$ is the directed graph with vertex set $S$ and edges $w \-> \lambda_q(w)$ for $w \in S$. The graph of action of  automata ${\A_1}_{q_1}, {\A_2}_{q_2}, \ldots, {\A_k}_{q_k}$ on $S$ is similarly defined as a directed graph with vertex set $S$ and edges $w \-> {\lambda_i}_{q_i}(w)$, $1 \leq i \leq k$ and $w \in S$.
\end{definition}

\begin{figure}[ht]
	\begin{subfigure}{0.499\textwidth}
		\centering	
		\includegraphics[width=0.99\textwidth]{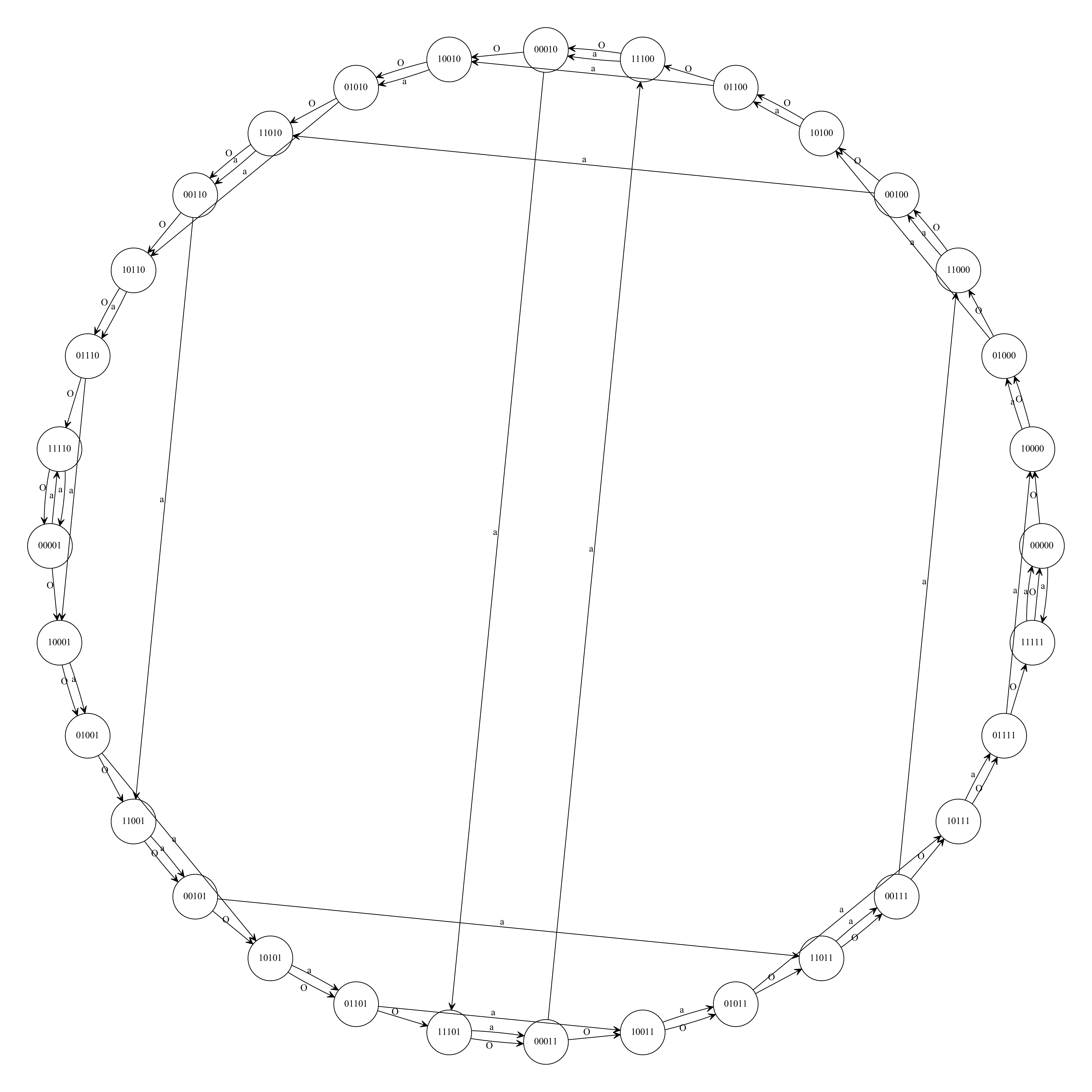}
		 \caption{The graph of action of adding machine and \zoran}
		 		\label{fig:zoran_act}
	\end{subfigure}
	\begin{subfigure}{0.499\textwidth}
		\centering	
		\includegraphics[width=0.99\textwidth]{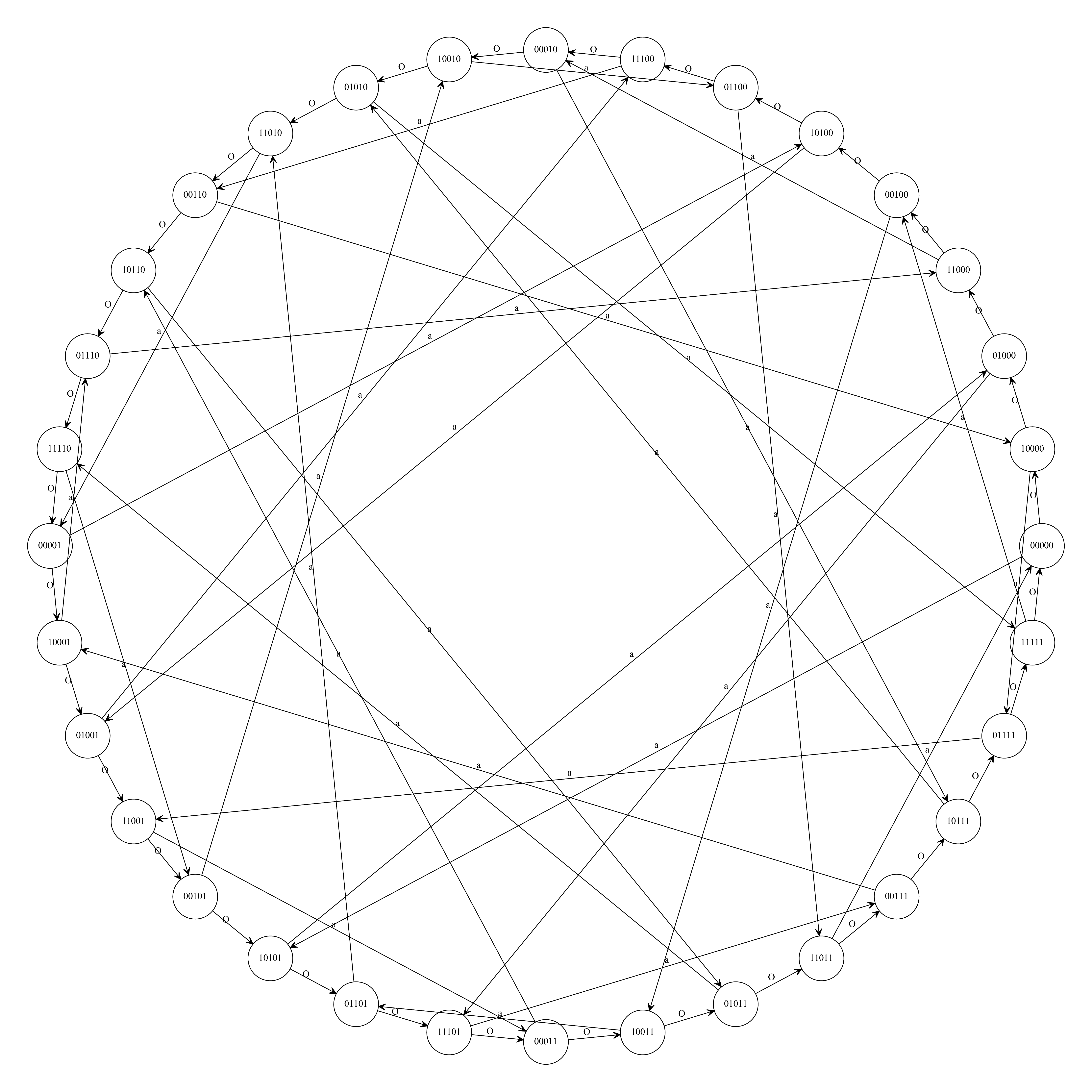}
		 \caption{The graph of action of adding machine and \alesh}
		  		 \label{fig:alesh_act}
	\end{subfigure}
	\caption{Examples of Schreier graphs}
	\label{fig:schreiergraphexamples}
\end{figure}

In this paper, we consider \bb{graphs of action on level $n$} of two automata, $\mathcal O$ and $\mathcal A$, with $\mathcal O$ being the adding machine (Figure \ref{fig:odom1}). Figure \ref{fig:schreiergraphexamples} shows examples of such graphs for $\mathcal A$ being \zoran{} (Figure \ref{fig:zoran}) and $\mathcal A$ being \alesh{} (Figure \ref{fig:alesh1}).

\begin{definition}\label{mooredef} A Mealy automaton is said to be a \bb{Moore machine} when the output does not depend on the last character of the input. That is, for all $s \in S$, $\lambda_s$ is constant: for all $x, y \in X$, $\lambda(s, x) = \lambda(s, y)$. In this case, we simply write $\lambda(s)$ for the value $\lambda_s$ takes.
\end{definition}
 
\begin{remark} 
In this definition, the output only depends on the current state $s$. Some authors use the definition of a Moore machine with a shift, where the output is determined by the ending state $\pi(s, x)$, and so does depend on the input.
\end{remark} 

Mealy automata $\mathcal A$ and $\mathcal B$ are said to be \bb{equivalent} if $\mathcal A(w) = \mathcal B(w)$ for all $w \in X^*$. 

\begin{definition} An initial Mealy automaton $\mathcal A$ is said to be \bb{minimal} if it has the smallest number of states among all the automata in its equivalence class. 
\end{definition}

Minimality is a classical notion, as is the algorithm that produces the minimal automaton in a given class; see \cite{sholomov} for a discussion of this algorithm (refer to \cite{GNS} for a discussion of this equivalence and a minimization algorithm in the more general case of asynchronous Mealy machines). 

Given automata $\mathcal A$ and $\mathcal B$ such that the output alphabet of $\mathcal A$ coincides with the input alphabet of $\mathcal B$, one can construct the \bb{product automaton}, denoted $\mathcal A \cdot \mathcal B$, which computes the composition $\mathcal A\circ \mathcal B$. We again refer to \cite{GNS} for the construction of the product automaton.

%

\subsection{Sections of tree endomorphisms}
 
\begin{definition} Let $g$ be an endomorphism of a $d$-regular rooted tree $\mathcal T$, and let $w$ be a finite word. A \bb{section} of $g$ by $w$, denoted $g|_w$, 
is an endomorphism $h$ of $\mathcal T$ such that for any word or sequence $v$, $g(wv) = g(w)h(v)$.
\end{definition}

\begin{remark}\label{statesectremark} a finite automaton $\mathcal A$ has only finitely many sections, which correspond to states in the connected component of the starting state in the diagram of the automaton $\mathcal A$.
More specifically, when $g$ is given by a Mealy machine $\mathcal A_q$, $g|_w = \mathcal A_{g(w)}$.
\end{remark}


With an invertible tree endomophism $g$ we can associate a \bb{portrait} diagram that uniquely determines $g$. For a finite word $w$, $g|_w$ acts on $X$ by a permutation when $g$ is invertible. The portrait consists of the infinite tree $\mathcal T$ with markings on the nodes: node corresponding to word $w$ is marked with the permutation of $X$ induced by $g|_w$. When $|X|=2$, we only mark nodes with the nontrivial permutation and leave others unmarked.

\begin{figure}
\centering
\includegraphics[width=0.25\textwidth]{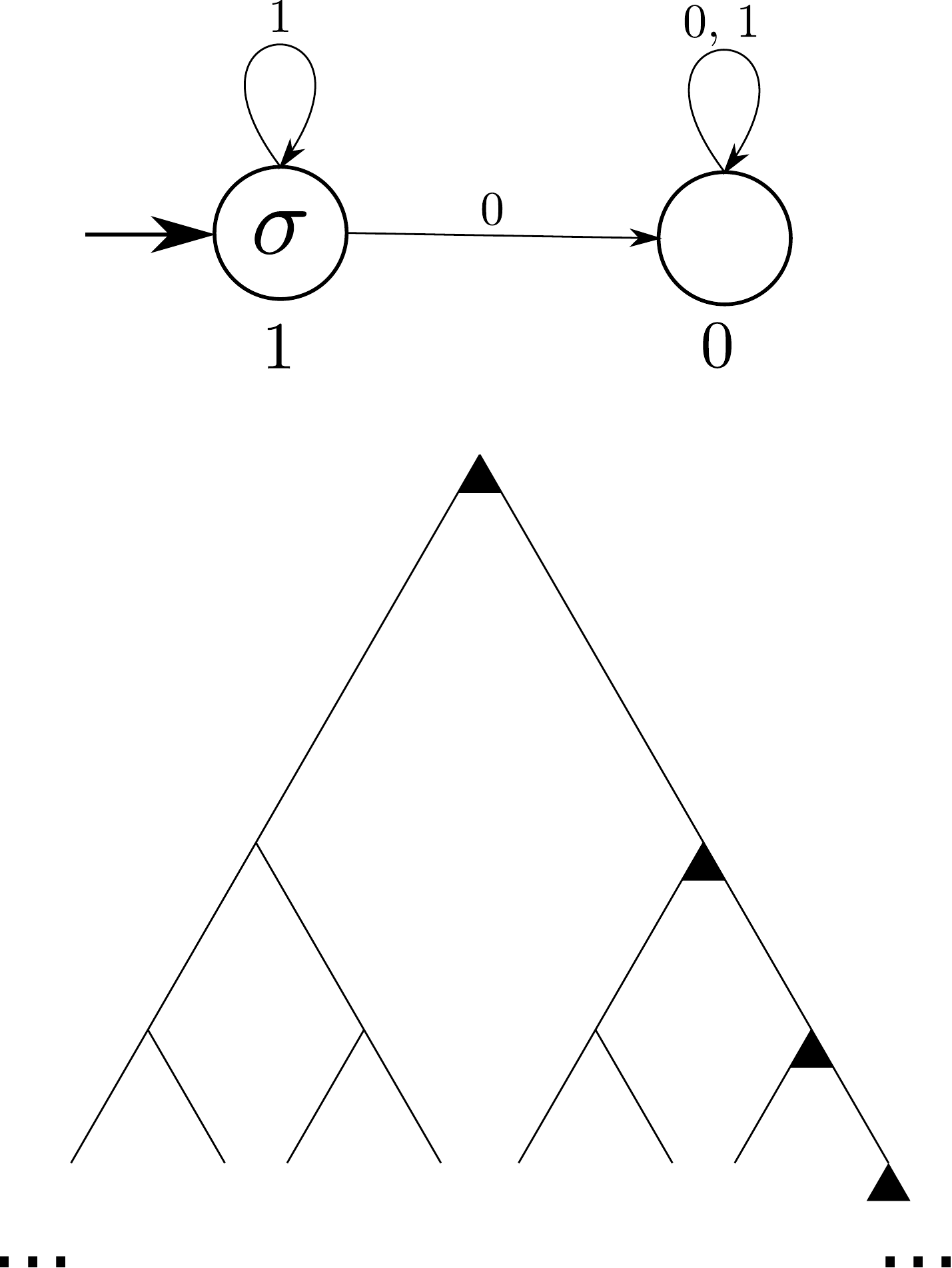}
\caption{The adding machine and its portrait}
\label{fig:odomportrait}
\end{figure}

\begin{example}
The portrait of the adding machine of Figure \ref{fig:odom1} is shown in Figure \ref{fig:odomportrait}. \tri
\end{example}

Every tree automorphism has a portrait, but not all tree automorphisms are given by finite automata. 
To any tree endomorphism $g$ we can associate a (possibly infinite) automaton $A=(S,g,X,\pi,\lambda)$ with the initial state labeled by $g$, such that the action of $A$ is identical to the action of $g$. We take $S=\{g|_w\ :\ w \in X^\nn\} \cup \{g\}$, and
define $\pi(h,x) \defeq h|_x$; $\lambda(h, x)=h(x)$.  This automaton of restrictions, in general, need not be finite. When it is finite, the tree automorphism $g$ is said to be \bb{finite-state}.

\begin{remark}
An automorphism $g$ of the tree $\mathcal T$ is finite-state if and only if its portrait contains a finite number of distinct (up to isomorphism of marked trees) subtrees. The subtrees in the portrait diagram define sections of $g$.
\end{remark}

We now prove several basic propositions related to sections of automorphisms which we use in subsequent chapters. These statements are well-known, but we include them for the reader's convenience.

\begin{proposition}\label{invsec}
If an endomorphism $g$ is invertible, then all of its sections are invertible, and for $w \in X^*$, $(g|_w)^{-1} = g^{-1}|_{g(w)}$.
\end{proposition}

\pf  Let $w \in X^*$ and $v \in X^\nn$. Then by definitions,
\begin{align*}
wv &=g^{-1}(g(wv))\\
&=g^{-1}(g(w)g|_w(v))\\
&=g^{-1}(g(w)) g^{-1}|_{g(w)}(g|_w(v))\\
&=w g^{-1}|_{g(w)}(g|_w(v)).
\end{align*}
Therefore  $g^{-1}|_{g(w)}(g|_w(v)) = v$, and the proposition holds. \sqr

\begin{proposition}\label{catsec}
Let $g$ be a tree endomorphism. Then for all $w, v$ finite words $w, v$ over $X$
 $g|_{wv} = (g|_w)|_v$.
\end{proposition}
\pf For any word $u$, 
\[
g(wvu) = g(w)g|_w(vu) =  g(w)g|_w(v)(g|_w)|_v(u) = g(wv)(g|_w)|_v(u).
\]
The proposition holds by definition. $\sqr$

\begin{proposition}\label{prosec}
Let $g$ and $h$ be tree endomorphisms. Then for all finite words $w$ over $X$,
 $(gh)|_w = g|_{h(w)}h|_w$.
\end{proposition}

\pf Let $v$ be a finite word. By the definition of section,
\begin{align*}
gh(wv) &= g(h(w)h|_w(v))\\
&= g (h(w)) g|_{h(w)}h_w(v),
\end{align*}
so $(gh)|_w = g|_{h(w)}h_w$. $\sqr$


\begin{corollary}\label{procatsec}
Let $g, h$ be tree endomorphisms. Then for any finite words $w$ and $v$,
 $(gh)|_{wv} = g|_{h(wv)}h|_{wv} = g|_{h(w)h|_w(v)}(h|_w)|_v$.
\end{corollary}

\begin{proposition}\label{powsec}
 $g^n|_w = g|_{g^{n-1}(w)} g|_{g^{n-2}(w)} \ldots g|_{g(w)} g|_w$.
\end{proposition}
\pf The result holds trivially when $n=1$. By Proposition \ref{prosec}, 
\begin{align*}
g^n|_w &= (g \circ g^{n-1})|_w\\
&= g|_{g^{n-1}(w)} (g^{n-1}|_w).
\end{align*}
The result follows by induction. $\sqr$

\begin{proposition}\label{gluesec}
Assume $g$ acts transitively on levels, $|w|=n$, and $a \in X$. 
Then $g^{2^n}|_w(a) \neq a$. 
\end{proposition}
\pf If $g^{2^n}|_w(a) = a$, then  $wa$, a word of length $n+1$, is a fixed point of $g^{2^n}$, contrary to the assumption that the length of the orbit of $g$ on words of length $n+1$ is $2^{n+1}$.

\subsection{Automata with bounded activity}

\begin{definition}\label{boundeddef} An automaton $A$ is said to have \bb{bounded activity} if the number of nontrivial sections on every level is bounded by a global constant $c$:
\[
\exists c: \forall n \in \nn:\  |\{ A|_w\ :\ A|_w \neq \bb{1}, w\in X^n\}| < c.
\]
\end{definition}

\begin{example}
The adding machine in Figure \ref{fig:odom1} has bounded activity. 

This automaton can also be defined by the portrait in Figure \ref{fig:odomportrait}, in which case it is clear that there is only one nontrivial section on every level.
$\tri$
\end{example}


The following proposition shows that the set of sections of powers of a bounded-activity automaton is finite if the powers are bounded by the number of words on the corresponding level. This fact is used to show Theorem \ref{dautom}.
\begin{proposition}\label{finprop}
If $A$ is a tree endomorphism given by a finite Mealy automaton $\mathcal{A}$ which is of bounded activity and acts transitively on levels, then the set 
\begin{align*}
\label{TAsectiondefn}
T_A := \{ A^n|_w\ :\ w \in X^*, n \leq 2^{|w|}\}
\end{align*}
is finite. 
\end{proposition}

\pf by Proposition \ref{powsec},
\[
T_A := \{ A|_{A^{n-1}(w)} A|_{A^{n-2}(w)} \ldots A|_{A(w)} A|_w\ :\ w \in X^*, n \leq 2^{|w|}\}.
\]

For a given $w$, consider a sequence of words $w, A(w), A^2(w), \ldots, A^{n-1}(w)$ with $n\leq 2^{|w|}$. By level transitivity of the action of $A$, all elements in it are distinct words of length $|w|$, and thus this sequence is a subset of vertices on level $|w|$.

Since $A$ is of bounded activity, there is a constant $c$ such that at most $c$ sections on every level are nontrivial. Hence the product
\[
A|_{A^{n-1}(w)} A|_{A^{n-2}(w)} \ldots A|_{A(w)} A|_w
\]
contains at most $c$ nontrivial factor. Since $A$ is finite-state by assumption, its nontrivial sections are enumerated by the finite set of states $S_A$ of $\mathcal{A}$.
Therefore, $|T_A| \leq |S_A|^c$. \sqr

\subsection{Measure-theoretic definitions}
We now give a few definitions relevant to probability theory and ergodic theory.

A \bb{cylinder set} $\cyl w$ is a clopen subset of $X^\nn$ given by
\[
\cyl w := \{wv\ :\ w \in X^*, v\in X^\nn\}.
\]

A \bb{probability vector} $p$ is a a vector $p : X \-> [0,1]$ with $\Sigma_{i \in X} p(i) = 1$. A \bb{stochastic matrix} on $X$ is a matrix $M: X \times X \-> [0,1]$ whose rows are probability vectors.

\begin{definition} \label{berndef}The \bb{Bernoulli measure} on $X^\nn$ defined by a probability vector $p$ is given on the cylinders $wX^\nn$ by
\[
\mu(wX^\nn) :=\prod_{i=0}^{|w|-1}p(w_i),
\]
and extended by the additivity properties on all Borel sets. The \bb{uniform Bernoulli measure} is given by $p=\left(\frac{1}{|X|},\ldots,\frac{1}{|X|}\right)$.
\end{definition}

Informally, this measures probability of a sequence of independent events (e.g. coin flips).

\begin{definition}\label{markdef} The \bb{Markov measure} defined by a probability vector $l=(l_x)$ of length $|X|$ and a stochastic matrix $L=(L_{x,y})$ of size $|X|\times |X|$ is given on the cylinder sets $\cyl w$ by
\[
\mu(\cyl w) :=l(w_0) \prod_{i=1}^{|w|-1}L_{w_{i-1},w_i}.
\]
\end{definition}

Informally, this measures the probability of events where the probability of an outcome may depend on what the preceding outcome was.

\subsection{Sections of a measure}

\begin{definition} The \bb{null measure} $\nu_0$ (or the \bb{trivial measure}) $\nu_0$ is the measure given by  $\nu_0(E) = 0$ for all measurable sets $E$.
\end{definition}

\begin{definition} Suppose $\mu$ is a probability measure on $X^\nn$.
If $\mu(wX^\nn) \neq 0$, then the \bb{section} of $\mu$ by the word $w \in X^*$, denoted $\mu|_w$, is the probability measure on $X^\nn$ uniquely defined by
\[
\mu|_w(\cyl v) := 
\displaystyle \frac{\mu(\cyl{wv})}{\mu(\cyl w)}
\]
for all $v \in X^*$. In the case $\mu(wX^\nn)=0$, we let $\mu|_w$ be the null measure.
\end{definition}

The section $\mu_w$ can be seen as the conditional probability given $w$.

For convenience, we also define sections for null measures: if $\mu$ is null, $\mu|_w = 0$ for all words $w$.

We say a word $w$ is \bb{admissible} (with respect to $\mu$) if the section of $\mu$ by $w$ is nontrivial (i.e. $\mu(wX^\nn) \neq 0$). We say a word $w$ is \bb{forbidden} if it is not contained in any admissible word.

Now we describe how to compute sections of measures.

\begin{proposition}\label{msec}
$\mu|_{wv} = (\mu|_w)|_v$ for all words $v, w \in X^*$.
\end{proposition}
\pf 

First, suppose $wv$ is not admissible, i.e. $\mu(wvX^\nn)=0$. Then either $w$ is also not admissible, or $w$ is admissible relative to $\mu$, but $v$ is not admissible relative to $\mu|_w$. Either way, $(\mu|_w)|_v$ is the null measure, and the proposition holds.

Now assume $\mu(wvX^\nn) \neq 0$; then $\mu(wvX^\nn) \neq 0$. For any word $u \in X^*$ we obtain
\begin{align*}
(\mu|_w)|_v(\cyl u) &= \frac{\mu|_w(\cyl{vu})}{\mu|_w(\cyl{v})}\\
&= \frac{\mu(\cyl{wvu})}{\mu(\cyl{w})\mu|_w(\cyl{v})}\\
&= \frac{\mu(\cyl{wvu})}{\mu(\cyl{wv})}\\
&= \mu|_{wv}(\cyl{u}). \sqr
\end{align*}

\begin{corollary}\label{sumsec}
Let $\mu = \displaystyle \sum_{i=1}^k a_i\mu_i$, where $a_i \geq 0$ and $\mu_i$ are probability measures. The for any admissible word $w$, 
\[
\mu|_w = \frac{1}{\mu \pcyl w}\sum_{i=1}^k a_i \mu_i \pcyl{w} \mu_i|_w.
\]
\end{corollary}
\pf For any word $v \in X^*$,
\begin{align*}
\mu|_w \pcyl v &= \frac{1}{\mu \pcyl w}
{\sum_{i=1}^k a_i\mu_i\pcyl{wv}}=\frac{1}{\mu \pcyl w}\sum_{i=1}^k a_i \mu_i\lp \cyl{w} \rp  \mu_i|_w \pcyl v.\ \ \sqr
\end{align*}


\section{Finite-state measures}
\begin{definition} A measure $\mu$ is \bb{finite-state} if admits only finitely many distinct sections.
\end{definition}

\begin{example}
Bernoulli and Markov measures (definitions \ref{berndef} and \ref{markdef}, respectively) are finite-state:
\begin{itemize}
\item any Bernoulli measure $\mu$ has only one (nontrivial) section: $\mu|_w = \mu$ whenever $w$ is admissible. Indeed,
let $p$ be the defining probability vector,
\begin{align*}
\mu|_w(\cyl v) &= \frac{\mu(\cyl{wv})}{\mu(\cyl w)}\\
 &= \frac{\prod_{i=0}^{|w|-1}p(w_i)\prod_{j=0}^{|v|-1}p(v_j)}{\prod_{i=0}^{|w|-1}p(w_i)}\\
 &= \prod_{j=0}^{|v|-1}p(v_j) = \mu(\cyl{v}).
\end{align*}

Note that if $p(x) = 0$ for some $x \in X$, then words $w$ containing $x$ are not admissible. Conversely, if $p$ is positive, then all words are admissible.

\item a Markov measure $\mu$ has at most $|X|+1$ nontrivial sections: $\mu$ (section by the empty word) and $\mu|_x$ for $x \in X$. This is because for all admissible words $w \in X^*$ and all $x \in X$,  $\mu|_{wx}=\mu|_x$. Indeed, assuming $w$ is not the empty word, we obtain
\begin{align*}
\mu|_{wa}(\cyl{v}) &= \frac{\mu(\cyl{wav})}{\mu(\cyl{wa})}\\
 &= \frac{\left(l(w_0)\prod_{i=1}^{|w|-1}L(w_{i-1},w_i)\right) L(w_{|w|-1},a)\left(L(a,v_0)\prod_{j=1}^{|v|-1}L(v_{j-1},v_j)\right)}{\left(l(w_0)\prod_{i=1}^{|w|-1}L(w_{i-1},w_i)\right)L(w_{|w|-1},a)}\\
 &= L(a,v_0)\prod_{j=1}^{|v|-1}L(v_{j-1},v_j)\\
 &= \frac{l(a)L(a,v_0)\prod_{j=1}^{|v|-1}L(v_{j-1},v_j)}{l(a)}\\
 &= \frac{\mu(\cyl{av})}{\mu(\cyl{a})}\\
 &=  \mu|_a(\cyl v).
\end{align*}
\end{itemize} \tri
\end{example}

\begin{definition} \label{kstepmarkovmeasuredefinition}A \bb{$k$-step Markov measure} is a measure $\mu$ such that for all words $v \in X^*$ of length $k$ and all words $w \in X^*$, $\mu|_{wv} = \mu|_{v}$ whenever $wv$ is admissible.
\end{definition}

Informally, this measures the probability of events where the probability of an outcome may depend on what the preceding $k$ outcomes was.

\begin{remark}
A Markov measure is a 1-step Markov measure. A $k$-step Markov measure on $X^\nn$ with $|X|=d$ is finite-state with at most $\displaystyle \frac{d^{k+1}-1}{d-1}$ sections.

Indeed, a finite $d$-tree of depth $k+1$ has $1 + d + d^2 + \ldots + d^k = \frac{d^{k+1}-1}{d-1}$ nodes, which encode all words of length not exceeding $k$. By definition, every nontrivial section of a $k$-step Markov measure is a section by one of these words. \sqr
\end{remark}

\begin{definition} \label{fsmeasdefin}
To any finite-state measure $\mu$ we associate an automaton $\mathcal A_\mu$ as follows.

Let $\mu_1, \ldots, \mu_n$ be the distinct sections of $\mu$.  Consider an automaton $\mathcal A_\mu$ with input alphabet $X$, output alphabet $Y\subset [0,1]$, state set $S=\{\mu_1, \ldots, \mu_n\}$, initial state $s_0 = \mu \in S$, and transition and output functions defined by
\begin{align}\label{measautomdef}
\pi(\mu_i, a) &:= \mu_i|_a;\\
\lambda(\mu_i, a) &:= \mu_i(\cyl a). \nonumber
\end{align}

We say that $\mathcal A_\mu$ \bb{determines} the measure $\mu$.
\end{definition}

\begin{proposition}\label{fsmeasautomprop}
The automaton $\mathcal A_\mu$ uniquely determines $\mu$ as follows: for any input word $w \in X^*$, the output word $\mathcal A_\mu(w)=p_0p_1\ldots p_{|w|-1}$ is a sequence of real numbers whose product is $\mu(\cyl w)$:
\begin{align}
\mu(\cyl w) = \prod_{i=0}^{|w|-1}(\mathcal A_\mu(w))_i. \label{fsproddef}
\end{align}
\end{proposition}
\pf The proposition holds for when $|w|=1$ by construction; assume it holds for all words of length $k$. Consider an arbitrary word  $w=w_0w_1\ldots w_k$ of length $k$. Then applying the inductive hypothesis, and then applying Proposition  \ref{msec} $k$ times, we obtain:
\begin{align*}
\prod_{i=0}^{k}p_i   &= \mu(\cyl{w_0w_1\ldots w_{k-1}}) \cdot ((\ldots(\mu|_{w_0})|_{w_1})|_{w_2})\ldots)|_{w_{k-1}}(\cyl{w_k})\\
&= \mu(\cyl{w_0w_1\ldots w_{k-1}}) \mu|_{w_0w_1\ldots w_{k-1}}(\cyl{w_k}) \\
&= \mu(\cyl{w_0w_1\ldots w_{k}})  = \mu(\cyl w),
\end{align*}
which completes the inductive step. \sqr

Note that if $\mu_i$ is a section of $\mu$, then the automaton defining $\mu_i$ can be obtained from $\mathcal A_\mu$ by changing the initial state to $\mu_i$ and possibly dropping some states (as some sectinos of $\mu$ might not be sections of $\mu_i$). 

\begin{definition} Suppose $\mu$ is a finite-state measure that admits a trivial section. We call the corresponding state of the defining automaton $\mathcal A_\mu$ trivial.
\end{definition}

\begin{figure}[ht]
\begin{subfigure}{\textwidth}
\centering
\includegraphics[width=0.5\textwidth]{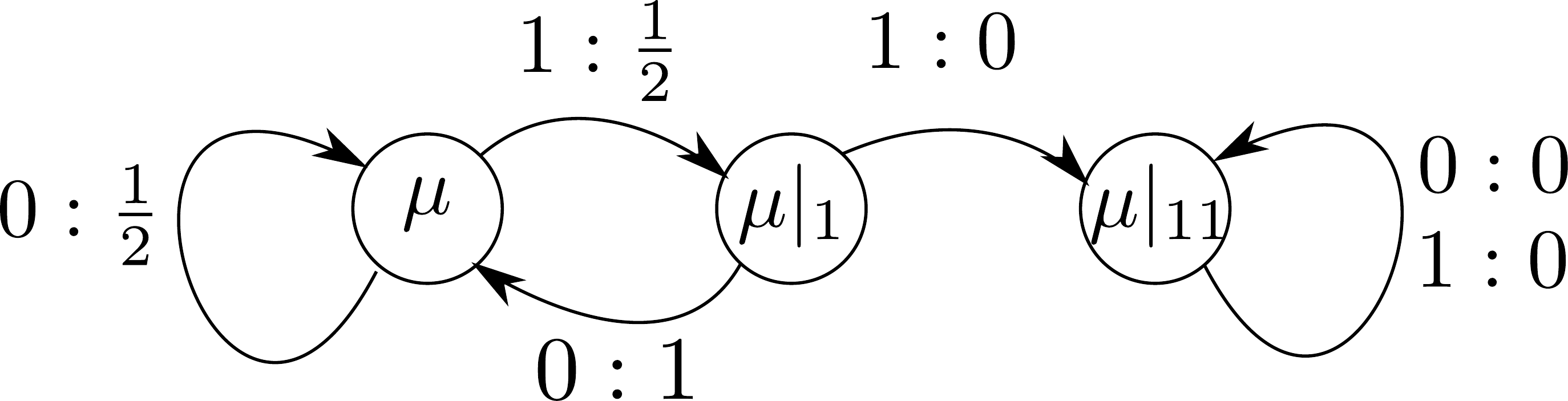}
\caption{2-step Markov measure $\mu$ on $\{0,1\}^\nn$}
\label{fig:fib_fs}
\end{subfigure}

\begin{subfigure}{\textwidth}
\centering
\includegraphics[width=0.7\textwidth]{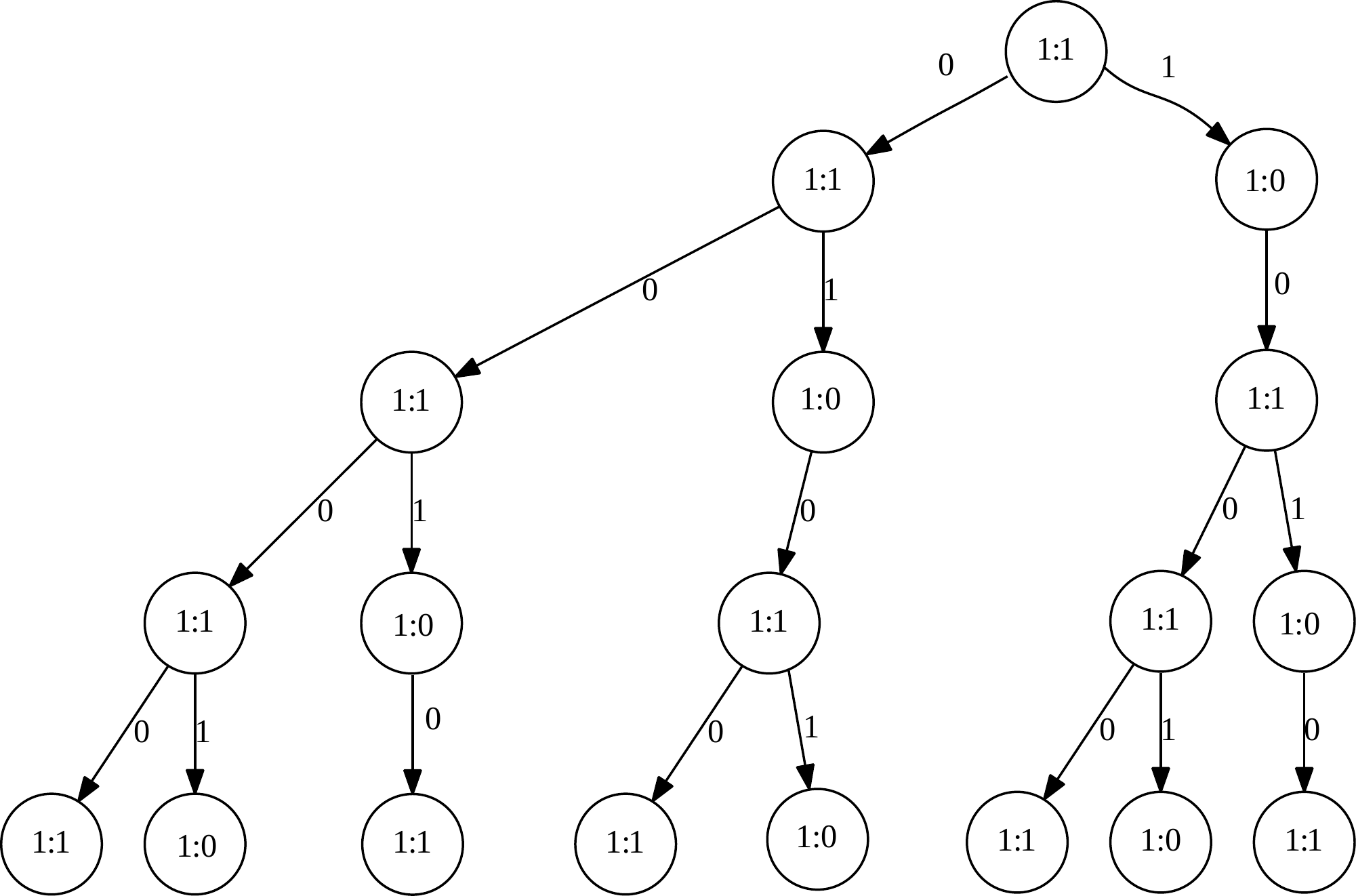}
\caption{Portrait of $\mu$ up to level $5$}
\label{fig:fib_portrait}
\end{subfigure}
\caption{A finite-state measure supported on the Fibonacci subshift}
\end{figure}

Refer to Example \ref{fibex} for an automaton with a trivial state; for example, the state $\mu|_{11}$ in Figure \ref{fig:fib_fs} is trivial.

\begin{remark}
Given a finite-state measure $\mu$, the automaton defined in \ref{measautomdef} is minimal and contains at most one trivial state.
\end{remark}

\begin{figure}[ht]
\begin{subfigure}{0.4\textwidth}
\centering
\includegraphics[width=0.6\textwidth]{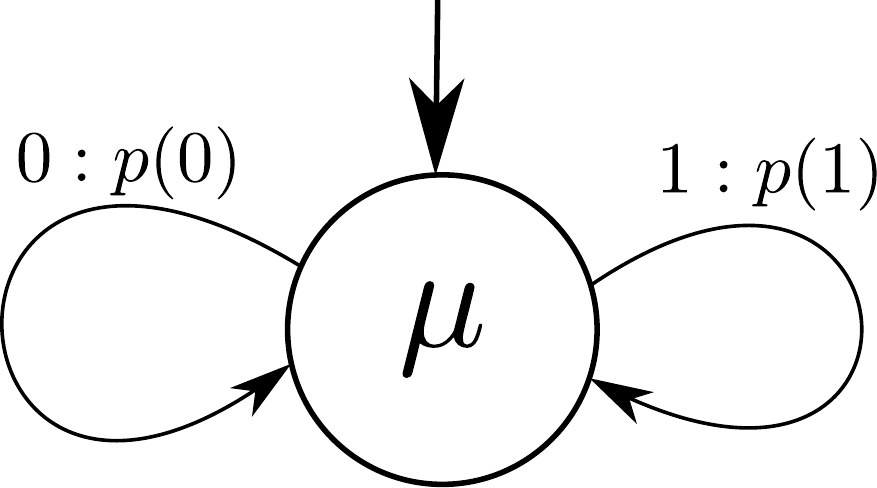}
\caption{Diagram of the automaton computing a Bernoulli measure}
\label{fig:bernoulli_fs}
\end{subfigure}
\begin{subfigure}{0.6\textwidth}
\centering
\includegraphics[width=0.8\textwidth]{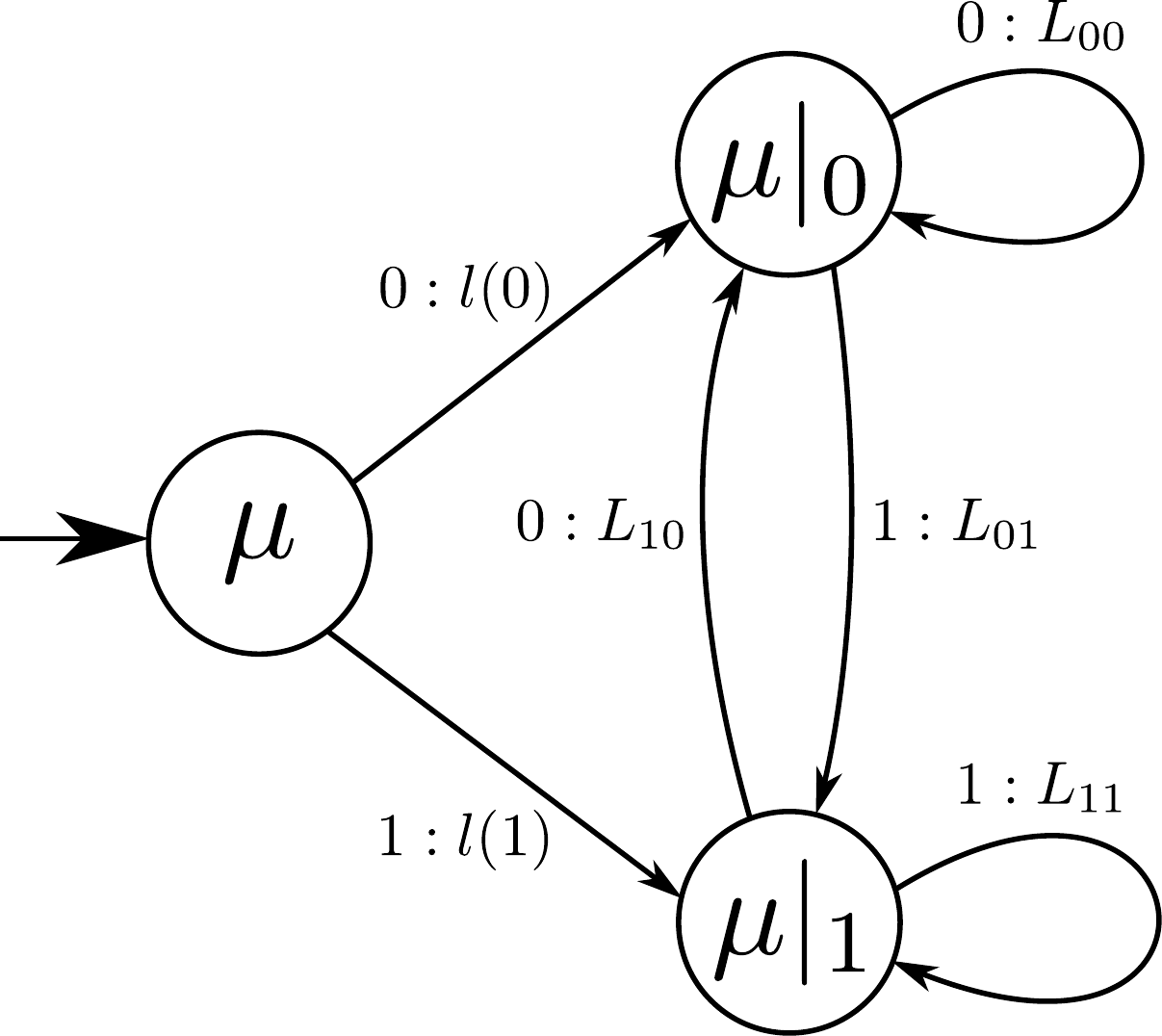}
\caption{Diagram of the automaton computing a Markov measure}
\label{fig:markov_fs2}
\end{subfigure}
\caption{Automata determining a Bernoulli and a Markov measure on $\{0,1\}^\nn$}
\end{figure}

\begin{example}
The automaton computing a Bernoulli measure on $\{0,1\}^\nn$ defined by a positive probability vector $p=(p(0),p(1))$ is depicted in Figure \ref{fig:bernoulli_fs}. \tri
\end{example}

\begin{example}
The automaton computing a Markov measure on $\{0,1\}^\nn$ defined by a positive probability vector $l=(l(0),l(1))$ and a positive probability matrix $L=(L_{ij})$ is depicted in Figure \ref{fig:markov_fs2}. \tri
\end{example}

\begin{figure}[ht]
\centering
\includegraphics[width=0.5\textwidth]{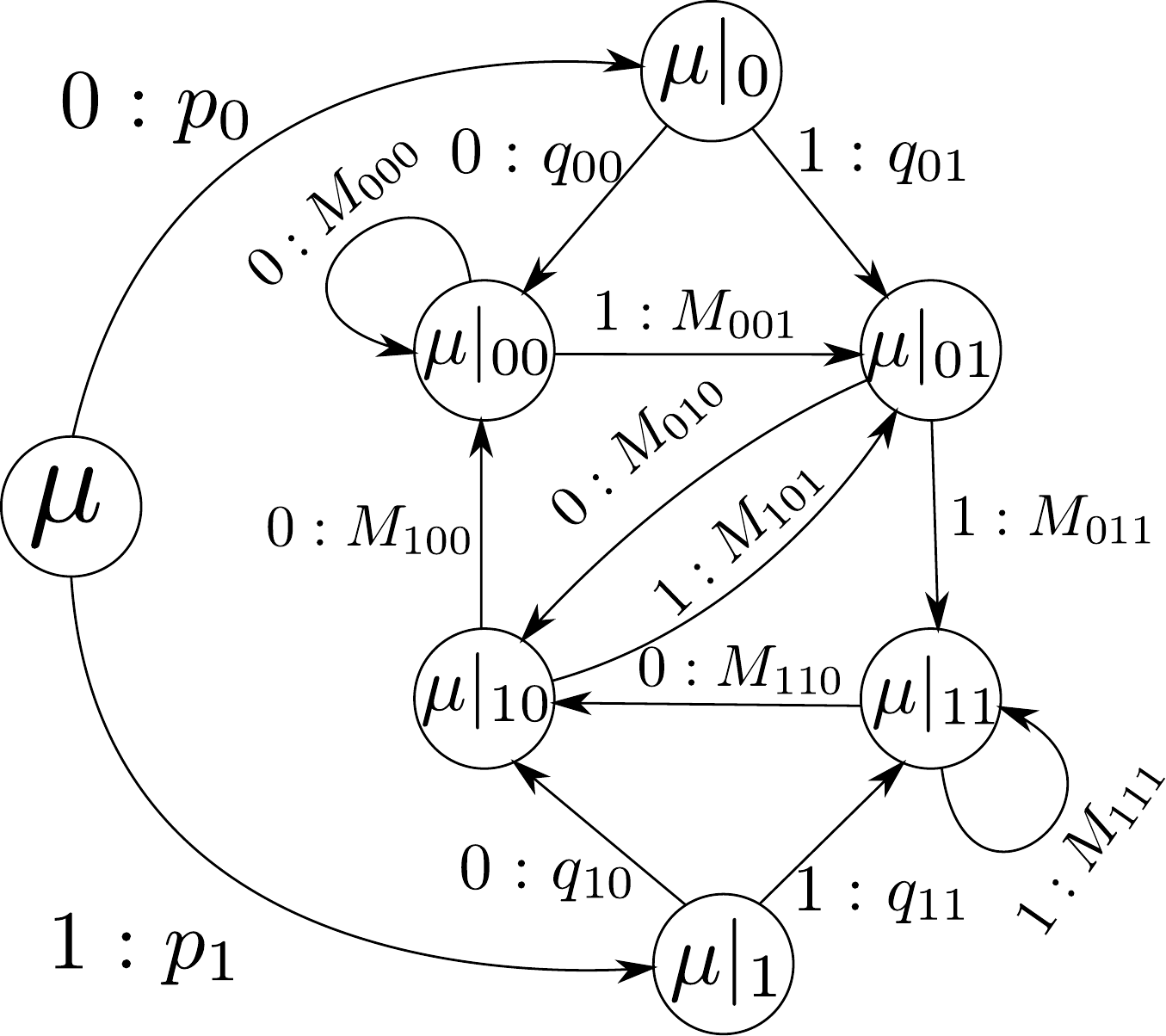}
\caption{Automaton defining a general $2$-step Markov measure on $\{0,1\}^\nn$}
\label{fig:2stepmarkov_fs}
\end{figure}

\begin{example}
Figure \ref{fig:2stepmarkov_fs} shows a general $2$-step Markov measure on $\{0,1\}^\nn$. Such a measure is determined by a probability vector $p$, a stochastic matrix $q$ and a probability tensor $M$ ($M_{ijk}$ gives the probability of $k$ given $ij$). \tri
\end{example}

Similarly to tree automorphisms, we define the \bb{portrait of the measure} $\mu$ to be the diagram consisting of the marked tree $\mathcal T$, where the node corresponding to a word $w$ is marked with the values $\mu|_w$ takes on cylinders $\cyl x$, $x\in X$.  A portrait defines a measure uniquely.

When dealing with probability measures, it is often convenient to consider the vector $p_w:=\left (\mu|_w(x_0X^\nn), \ldots, \mu|_w(x_{d-1}X^\nn \right)$ up to scaling. Since $\sum_{i=0}^{d-1} \mu(x_iX^\nn)=1$,
the \bb{proportion} $p_{w_0} : p_{w_1} : \ldots : p_{w_{n-1}} \in \rr P^{d-1}$ determines the values of $\mu$ on $X$ unambiguously. We then use the  proportion as the corresponding label in the portrait. 

\begin{figure}[ht]
\centering
\includegraphics[width=0.75\textwidth]{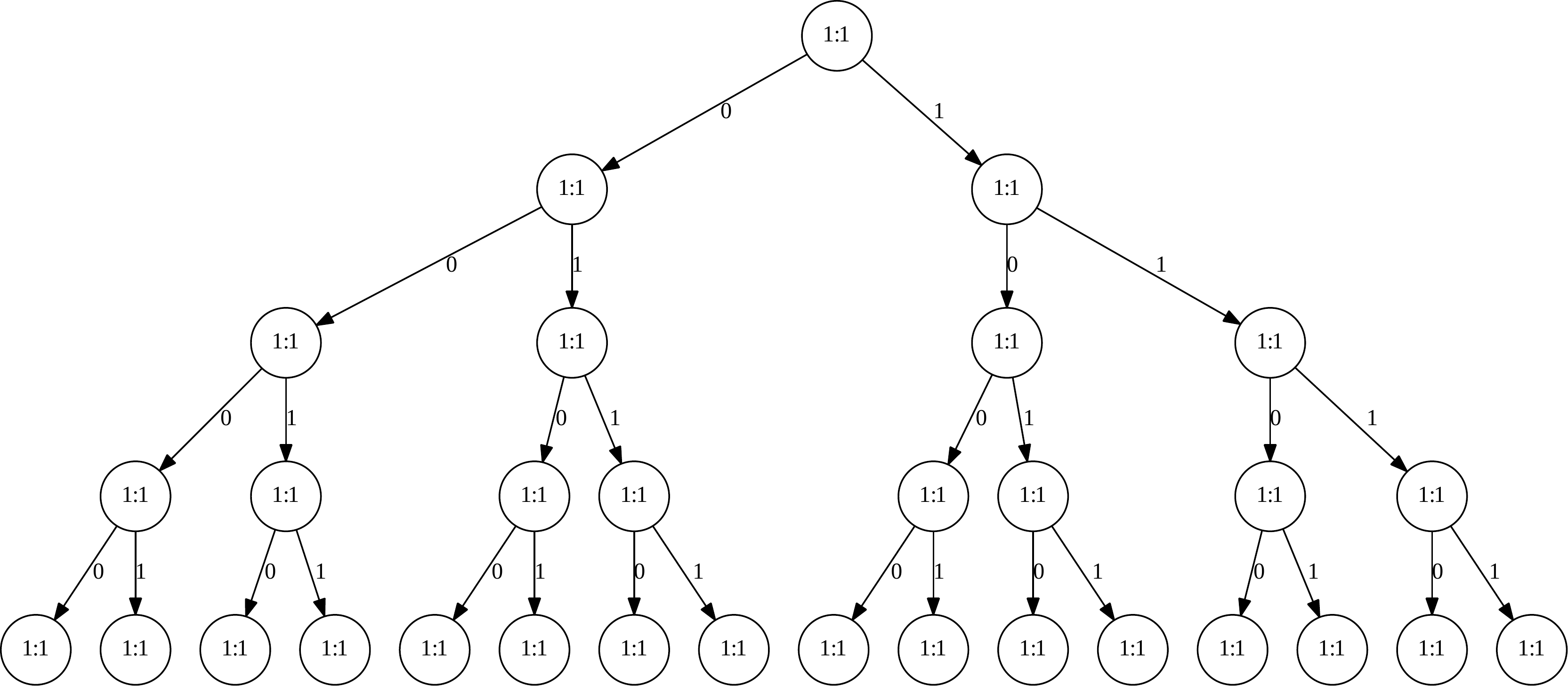}
\caption{Portrait of the uniform Bernoulli measure up to level $5$}
\label{fig:uniform_benoulli_portrait}
\end{figure}

\begin{example}
The uniform Bernoulli measure on a binary alphabet has one section whose proportion is $1:1$. Its portrait is shown in Figure \ref{fig:uniform_benoulli_portrait}. \tri
\end{example}

\begin{remark}
As with automorphisms, one can draw the portrait of any probability measure on the space $X^\nn$, but not all probability measures are finite-state.

Again, as with automorphisms, even small automata define interesting finite-state measures.
\end{remark}

It should be noted that even small automata define interesting finite-state measures.

\begin{example}\label{fibex}
The measure $\mu$ defined by the automaton in Figure \ref{fig:fib_fs} is a $2$-step Markov measure on $\Omega =\{0,1\}^\nn$ that
is not a $1$-step Markov measure on $\Omega$. It is supported on the \bb{Fibonacci subshift}, which is the (shift-invariant) subset of $\Omega$ consisting of all sequences that do not contain consecutive $1$'s. The number of nontrivial sections of $\mu$ by words of length $n-1$ is the $n$'th Fibonacci number as can be seen in the portrait of $\mu$ shown in Figure \ref{fig:fib_portrait}
(to simplify the figure, we omitted the subtrees corresponding to the null measure). \tri
\end{example}

\section{Images of finite-state measures under tree automorphisms}\label{imfsmeas_thesissection}
Given a finite-state measure $\mu$ and a tree automorphism $g$, we consider the pushforward measure $g_*\mu$ defined by $g_*\mu(E) = \mu(g^{-1}(E))$ for all measurable sets $E$. We say that $g_*\mu$ is the image of $\mu$ under (the map) $g$.

The following proposition is useful for constructing the automata of finite-state measures which are images under automaton automotphisms. 

\begin{proposition}\label{imsec}
Let $A=(X,S,s_0,\pi,\lambda)$ be a Mealy automaton with initial state $s_0=g$ acting on $\mathcal T$, and let $\nu$ be a probability measure on $\partial \mathcal T$. Then for $x \in X$,
\begin{align*}
(g_*\nu)(\cyl x) &= \sum_{y \in \lambda_g^{-1}(x)} \nu(\cyl y);\\
(g_*\nu)|_x &= \frac{\displaystyle\sum_{y \in \lambda_g^{-1}(x)} \nu(\cyl y)\pi_g(y)_*(\nu|_y)}
{\displaystyle\sum_{y \in \lambda_g^{-1}(x)} \nu(\cyl y)}.
\end{align*}
where $\pi_g(x) := \pi(g, x)$ and $\lambda_g(x) := \lambda(g,x)$.
\end{proposition}
\pf Note that for a word $w \in X^*$, 
\[
g^{-1}(\cyl{xw}) = \bigsqcup_{y \in \lambda_g^{-1}(x)} y \pi_g(y)^{-1}(\cyl{w}).
\]
By definition, 
\begin{align*}
(g_*\nu)|_x(\cyl w) &= \frac{(g_*\nu)(\cyl{xw}}{(g_*\nu)(\cyl{x})}\\
&= \frac{\nu(g^{-1}(\cyl{xw})}{\nu(g^{-1}(\cyl{x}))}\\
&= \frac{\displaystyle\sum_{y \in \lambda_g^{-1}(x)} \nu(\cyl y) \nu|_y \left(\pi_g(y)^{-1}(\cyl{w})\right)}
{\displaystyle\sum_{y \in \lambda_g^{-1}(x)} \nu(\cyl y) \nu|_y \left(\pi_g(y)^{-1}(\cyl{})\right)}\\
 &= \frac{\displaystyle\sum_{y \in \lambda_g^{-1}(x)} \nu(\cyl y)\pi_g(y)_*(\nu|_y)(\cyl{w})}
{\displaystyle\sum_{y \in \lambda_g^{-1}(x)} \nu(\cyl y)}. \sqr
\end{align*}
\begin{corollary}\label{imsecber}
When $g$ is as in Prop. \ref{imsec}, and $\nu$ is a Bernoulli measure given by probability vector $p$, then its image under $g$ satisfies
\[
(g_*\nu)|_x = \frac{\displaystyle\sum_{y \in \lambda_g^{-1}(x)} p(y)\pi_g(y)_*(\nu)}
{\displaystyle\sum_{y \in \lambda_g^{-1}(x)} p(y)}.
\]
In particular, when $\nu$ is uniform Bernoulli, $(g_*\nu)(\cyl x) = |\lambda_g^{-1}(x)|/|X|$, and
\begin{align*}
(g_*\nu)|_x &= \frac{1}{|\lambda_g^{-1}(x)|}\sum_{y \in \lambda_g^{-1}(x)}{\pi_g(y)_*(\nu)}.
\end{align*}
\end{corollary}
When $\nu$ is uniform Bernoulli, its pushforwards by invertible endomorphisms are also uniform Bernoulli:

\begin{proposition}\label{invpush}
When $\nu$ is uniform Bernoulli and $g$ is invertible, $g_*\nu = \nu$.
\end{proposition}
\pf For $w \in X^*$, 
\begin{align*}
g_*\nu(\cyl{w}) &= \nu(g^{-1}(\cyl{w}) = \nu(\cyl{g^{-1}(w)}) = |X|^{-|w|}=\nu(\cyl{w}). \sqr
\end{align*}


\section{Log map}

Let $A$ be an automorphism of the $d$-regular rooted tree $\mathcal T$ that acts transitively on each level. 
Recall that level $n$ of the tree consists of $d^n$ words of length $n$.
Hence for any pair of words $w_1, w_2$ of length $n$, there is a unique integer $k$, $0 \leq k \leq d^{n}-1$ such that $A^k(w_1) = w_2$. Furthermore, if $A^{k}(w_1) = A^{k'}(w_1)$ for some integers $k$ and $k'$, then $k\equiv k' \mod d^n$. 

\begin{definition} For any $n \geq 1$, the \bb{displacement} function $\bb{d}_{A,n} : X^n \times X^n \-> \zz/d^n\zz$ is defined on pairs of words $w_1, w_2$ of length $n$ by
\[
\bb{d}_{A,n}(w_1, w_2) := [k]_{d^n},
\]
where $A^k(w_1)=w_2 $ and $[k]_{d^n} \in \zz/d^n \zz$ is the equivalence class mod $d^n$. We write $[k]$ when $n$ is understood from the context.
\end{definition}

\begin{definition}For any integers $m$ and $n$, $1 \leq m \leq n$, the \bb{natural projection} $\phi_{m,n}:\zz/d^n\zz \-> \zz/d^m\zz$ is defined by $\phi_{m,n}([k]_{d^n}) := [k]_{d^m}$. These functions are homomorphisms of rings.
\end{definition}

The functions $\bb{d}_{A,n}$ for different values of $n$ are compatible with each other with respect to the natural projections.

\begin{proposition}\label{dcompat}
Suppose $|w_1|=|w_2|=n$ and $a, b \in X$. Then 
\[
\phi_{n,n+1}(\bb{d}_{A,n+1}(w_1a, w_2b)) = \bb{d}_{A,n}(w_1, w_2).
\]
\end{proposition}
\bb{Proof:} Let $\bb{d}_{A,n}(w_1, w_2) = [k]$  so that $A^k(w_1)=w_2$, with $0 \leq k \leq d^n-1$. 
Let $a'=A^k|_{w_1}(a)$. Then $A^k(w_1a) = w_2a'$. Note that
\begin{align*}
A^{d^{n+k}}(w_1a) &= A^{d^n}(w_2a')\\
&= A^{d^n}(w_2)A^{d^n}|_w(a').
\end{align*}

By Proposition \ref{gluesec},
\[
A^{d^n}|_w(a'), A^{2d^n}|_w(a'), \ldots, A^{(k-1)d^n}|_w(a'), 
\]
are all distinct. Since $|X|=d$, this implies $A^{td^n}|_w(a')=b$ for some $t$, $0 \leq t \leq d-1$. Thus
$A^{td^n+k}(w_1a) = w_2b$, whence $d_{A,n+1}(w_1a, w_2b)=[k+td^n]$.

Since $\phi_{n, n+1}([k+td^n]) = [k]$, the proposition holds. $\sqr$.

In addition to the tree endomorphism $A$, let us consider a tree endomorphism $B$.

\begin{definition} \label{logAndefn}  For any $n \geq 1$, 
 $\Log_{A,n}(B) : X^n \-> \zz/d^n\zz$ is a function which calculates the displacement of a word $w$ of length $n$ along the orbit of $A$ under the action of $B$:
\begin{align*}
\Log_{A,n}(B)(w) &:= \bb{d}_{A,n}(w,B(w)).
\end{align*}
\end{definition}

Note that for any word $w$ of length $n$, $A^{\Log _{A,n}(B)}(w) = B(w)$, which motivates the name ``logarithm'' for this function.

 \begin{corollary}\label{logcorollary}  For any word $w$ of length $n$ and character $x \in X$,
\[
\phi_{n,n+1}(\Log_{A,n+1}(B)(wa)) = \Log_{A,n}(B)(w).
\]
\end{corollary}
 
In other words, the displacement of $wa$ by $B$ along the orbit of $A$ is either the same as displacement of $w$ or differs from it by a multiple of $d^n$.

\begin{corollary}\label{comdiag1} For any integers $m$ and $n$, $1 \leq m < n$, the following diagram commutes:
\[
\begin{diagram} 
X^n & \rTo^{\sigma_r^{n-m}} & X^m  \\
 \dTo^{\Log_{A,n}(B)}&  &\dTo^{\Log_{A,m}(B)} \\
\zz/d^n\zz &  \rTo^{\phi_{m,n}} & \zz/d^m\zz \\
 \end{diagram}
\]
Here, $\sigma_r$ is the operator that trims the word, deleting the least letter: $\sigma(wa)=w$ for any word $w$ and a character $a\in X$.
\end{corollary}
\pf This follows from Corollary \ref{logcorollary} by induction on $n-m$. \sqr

Let $\zz_d$ be the inverse limit of the directed system 
\[
\begin{diagram}
\zz/d^n\zz &  \rTo^{\phi_{m,n}} & \zz/d^m\zz \\
 \end{diagram}
\]
(for $m,n \in \nn$). $\zz_d$ comes with a natural structure of a ring, and is known as the ring of the $d$-adic integers (note that $d$ need not be prime). 

Since the boundary of the tree $\partial \mathcal T$ can also be seen as the inverse limit of the directed system
\[
\begin{diagram} 
X^n & \rTo^{\sigma_r^{n-m}} & X^m,
 \end{diagram}
\]
Corollary \ref{comdiag1} implies that 
there exists a unique function $\Log_A(B) : \partial \mathcal T \-> \zz_d$, which restricts to $\Log_{A,n}(B)$ on level $n$ for all $n$.

\begin{definition} \label{logdefinition}
The logarithm $\Log_A(B)$ is the inverse limit
\[
\Log_A(B) = \varprojlim_n \Log_{A,n}(B).
\]
That is, it is the unique function $\Log_A(B):\partial \mathcal T \-> \zz_d$ that makes the following diagram commute:
\[
\begin{diagram} 
\partial \mathcal T & \rTo^{\pi_n} & X^n  \\
 \dTo^{\Log_{A}(B)}&  &\dTo^{\Log_{A,n}(B)} \\
\zz_d&  \rTo^{\pi_n} & \zz/d^n\zz \\
 \end{diagram}
\]
($\pi_n$ are the natural projections of the corresponding inverse limits). 
\end{definition}

Any positive integer $N$ admits a unique $d$-ary expansion
\[
N = \sum_{i=0}^k a_i d^i,
\]
where each $0 \leq a_i \leq d-1$. This way, the set $\zz/d^n\zz$ can be identified with the set of words of length $n$ over the alphabet $X$. Consequently, the set $\zz_d$ can be identified with infinite words in alphabet $X$, which, in term, are identified with the boundary of the tree $\partial \mathcal T$. Therefore, $\Log_{A,n}(B)$ can be seen as a transformation of the $n$'th level, and $\Log_A(B)$ can be regarded as a transformation of $\partial \mathcal T$. 

\begin{proposition} \label{logendomprop}
There exists an endomorphism of the tree $\mathcal T$ such that $\Log_{A,n}(B)$ is the restriction of the endomorphism to level $n$, and $\Log_{A}(B)$ is the action of the endomorphism on the boundary $\partial \mathcal T$. 
\end{proposition}

\pf Let $L : X^* \-> X^*$ be the transformation that coincides with $\Log_{A,n}(B)$ on level $n$ for all $n$. By construction, $L$ preserves the levels.
Corollary \ref{logcorollary} implies that $L$ maps adjacent vertices of $\mathcal T$ to adjacent vertices. Therefore, $L$ is an endomorphism. 
By Definition \ref{logdefinition}, the action of $L$ on the boundary $\partial \mathcal T$ is exactly $\Log_A(B)$. \sqr


In the rest of the paper, we deal with $d=2$, and so identify (and use interchangeably) the dyadic numbers and infinite binary sequences (elements of $\partial \mathcal T$).
 
\begin{remark}
The construction of the logarithm map $\Log_A(B)$ (including Proposition \ref{logendomprop}) can be extended from $d$-regular trees to spherically homogeneous trees (defined in, e.g., \cite{BORT}).
\end{remark} 

\section{The automaton computing the Log map}
Here and onwards we assume that $X=\{0,1\}$, that is, $\mathcal T$ is the binary rooted tree.

Let $A$ be an automorphism of the tree that acts transitively on each level, and let $B$ be an endomorphism. In light of Proposition \ref{logendomprop}, the Log map $\Log_{A}(B)$ can be regarded as an endomorphism of $\mathcal T$.

To simplify notation, we will denote by $\Log$ both $\Log_A(B)$ and $\Log_{A,n}(B)$.

In this section we construct an automaton which computes this endomorphism.

We further assume that the automorphism $A$ is of bounded activity (in the sense of Definition \ref{boundeddef}). An example of such endomorphism is the adding machine, whose automaton is shown in Figure \ref{fig:odom1}.

\begin{remark}
Any tree automorphism that acts transitively on levels is conjugate to the adding machine.
\end{remark}

The assumption that $A$ is of bounded activity allows us to prove the following lemma, which will be useful in the construction of the automaton for  $\Log$.

\begin{lemma}\label{fincor}
If $A$, $B$ are tree endomorphisms given by finite automata, $A$ is  of bounded activity and acts transitively on all levels, then the set $S_{A,B}$  consisting of triples of sections:
\[
S_{A,B} := \{ (B|_w, A^{d(w)}|_w, A^{2^{|w|}}|_w)\ :\ w \in X^* \}
\]
is finite.
\end{lemma}
\pf By Proposition \ref{finprop}, the set 
\[
T_A := \{ A^n|_w\ :\ w \in X^*, n \leq 2^{|w|}\}
\]
is finite. Note that $A^{d(w)}|_w, A^{2^{|w|}}|_w \in T_A$ for all $w \in X^*$. Let $S_B$ be the set of states of the automaton of $B$. Then
\[
 |S_{A,B}| \leq |S_B| \cdot |T_A|^2.  \sqr
\]

The set $S_{A,B}$ is going to be the set of states of our automaton. See Example \ref{aleshdex} for an explicit computation of $S_{A,B}$.

\begin{theorem}\label{dautom}
Let $A$, $B$ be as above. Consider the automaton $\L=\L_{A,B}$ with set of states $S_{A,B}$, initial state $(B, \bb{1}, A)$ (where $\bb{1}$ is the identity automorphism), 
transition function $\pi$ defined by 
\begin{align*}
\pi((\beta, \gamma, \delta), a) &:=  (\beta', \gamma', \delta'), \mbox{where} \\
\beta' &= \beta|_a\\
\gamma' &= \begin{cases}
\gamma|_a  \mbox{ if } \beta(a) = \gamma(a);\\
(\gamma\delta)|_a, \mbox{ otherwise};  
\end{cases}\\
\delta' &= \delta^2|_a\\
\end{align*}
and the output function  $\lambda$ given as follows:
\begin{align*}
 \lambda((\beta, \gamma, \delta), a) &:=
 \begin{cases}
  0, \mbox{ if } \beta(a) = \gamma(a);\\
  1, \mbox{ otherwise}.  
 \end{cases}
\end{align*}
Then the transition function is well-defined, and the automaton $\L$ outputs $\Log_{A,n}(B)$ as a dyadic integer:
\begin{align}
\label{logautomoutput}
\Log_{A,n}(B)(w) = \sum_{i=0}^{|w|-1}\L(w)_i 2^i,
\end{align}
where $n=|w|$, $\Log_{A,n}(B)$ is the displacement function in Definition \ref{logAndefn}, and $\L(w)_i$ is the $i$'th character of the word $\L(w)$.
\end{theorem}

\pf We first show that upon reading a word $w$, the automaton $L$ ends up in the state
\[
 \left(B|_w, A^{d(w)}|_w, A^{2^{|w|}}|_w  \right) \in S_{A,B}.
\]

This hypothesis holds for the empty word. We proceed by induction on $|w|$. 

Assume the hypothesis holds for all $|w| \leq  n$. 

To prove the inductive hypothesis for words of length $n+1$, let $|w|=n$ and $a \in X$, and assume that 
$\L$ is in the state $(\beta, \gamma, \delta)$ after reading  $w$. We show that
\[
(\beta',\gamma',\delta'):= \pi((\beta, \gamma, \delta), a) = \left(B|_{wa}, A^{d(wa)}|_{wa}, A^{2^{|wa|}}|_{wa}  \right).
\]

Indeed:
\begin{enumerate}
 \item $\beta' = \beta_a$ by definition, and 
 \begin{align*}
   B|_{wa} &= (B|_w)|_a  \mbox{\ \ (by Proposition \ref{catsec})}\\
   &=\beta_a\\
   &=\beta'.
 \end{align*}
\item Note that $A^{2^{|w|}}(w)=w$ by transitivity of $A$. By definition, $\delta' = \delta^2|_a = \delta|_{\delta(a)}\delta|_a$. Now
  \begin{align*}
   A^{2^{|wa|}}|_{wa} &= A^{2^{|w|+1}}|_{wa}\\
   &= (A^{2^{|w|}})^2|_{wa}\\
   &= A^{2^{|w|}}|_{A^{2^{|w|}}(wa)}A^{2^{|w|}}|_{wa} \mbox{\ \ (by Proposition \ref{prosec})}\\
   &= A^{2^{|w|}}|_{A^{2^{|w|}}(w)A^{2^{|w|}}|_w(a)} (A^{2^{|w|}}|_w)|_a\\
   &= A^{2^{|w|}}|_{w \delta(a)} \delta|_a \mbox{\ \  (by Proposition \ref{catsec}, inductive assumption, and $A^{2^{|w|}}=w$)}\\
   &= (A^{2^{|w|}}|_w)|_{\delta(a)} \delta|_a \\
   &= \delta|_{\delta(a)}\delta|_a\\
   &=\delta'.
  \end{align*}
\item Let $d(w) := \Log_{A,n}(B)(w)$. By definition of $\Log_{A,n}(B)$, $B(w) = A^{d(w)}(w)$. Note that
\begin{align*}
   B(wa) &= B(w)B|_w(a) = A^{d(w)}(w)\beta(a);\\
   A^{d(w)}(wa) &= A^{d(w)}(w)A^{d(w)}|_w(a) = A^{d(w)}(w)\gamma(a).
\end{align*}
If $\beta(a)=\gamma(a)$, then $B(wa) = A^{d(w)}(wa)$, and thus $d(wa)=d(wa)=d(w)$ by definition of $d = \Log_{A,n}(B)$.
Otherwise, $d(wa)=d(w)+2^{|w|}$ since this is the only other possibility.
Therefore,
\[
 A^{d(wa)}=\begin{cases}
            A^{d(w)}, \mbox{ if } \beta(a)=\gamma(a);\\
            A^{d(w)}A^{2^{|w|}}, \mbox{ otherwise}.
           \end{cases}
\]
Now we compute:
\begin{align*}
 A^{d(w)}|_{wa} &= (A^{d(w)}|_w)|_a\\
&= \delta|_a;\\
 A^{d(w)}A^{2^{|w|}}|_{wa} &= A^{d(w)}|_{A^{2^{|w|}}(wa)} (A^{2^{|w|}}|_w)|_a\\
&= A^{d(w)}|_{A^{2^{|w|}}(w) A^{2^{|w|}}|_w(a)} \delta|_a\\
&= A^{d(w)}|_{w \delta(a)} \delta|_a\\
&= (A^{d(w)}|_w)|_{\delta(a)} \delta|_a\\
&= \gamma|_{\delta(a)} \delta|_a\\
&= (\gamma\delta)|_a
\end{align*}

Therefore
\[
A^{d(wa)}|_{wa} = 
\begin{cases}
\gamma|_a  \mbox{ if } \beta(a) = \gamma(a);\\
(\gamma\delta)|_a \mbox{ otherwise}.  
\end{cases}
\]
This matches the definition of $\gamma'$, and thus $\gamma'=A^{d(wa)}|_{wa}$.
\end{enumerate} 

In particular, we have verified that the transition function $\pi$ is well-defined, since its values are always in the set $S_{A,B}$.

This completes the proof of the hypothesis that the automaton is in state $\left(B|_w, A^{d(w)}|_w, A^{2^{|w|}}|_w  \right)$ after reading $w$.

Furthermore, we observed that
\[
d(wa)= \begin{cases}
  d(w), \mbox{ if } \beta(a) = \gamma(a);\\
  d(w) + 2^{|w|} \mbox{ otherwise}.  
 \end{cases}\\ 
\]
From this observation and the definition of $\lambda$, equation \ref{logautomoutput} follows by induction.

This completes the proof of the theorem. $\sqr$

\begin{proposition}\label{dmoore}
When  $A$ and $B$ are as in Theorem \ref{dautom} and, additionally, $B$ is invertible, the automaton $L_{A,B}$ is a Moore machine (as in Definition \ref{mooredef}). Recall that the value of the output function $\lambda(s,x)$ of a Moore machine only depends on the state $s$.
\end{proposition}

\pf By assumption, $A$ is invertible, and so is $A^{d(w)}$ for any $w \in X^*$. $B$ is invertible by assumption.
By Proposition \ref{invsec}, their sections $\beta = B|_w$ and $\gamma = A^{d(w)}|_w$ are invertible, and so is $\beta\gamma^{-1}$. 

Now the set of permutations $\op{Perm}(\{0,1\}) =\{\bb{1}, \sigma\}$, so either $\beta\gamma^{-1}(x) = (x)$, or $\beta\gamma^{-1}(x) = \sigma(x)$.

In the first case, $\lambda(\beta, \gamma, \delta)(x) = 0$ for $x\in \{0,1\}$. 

Otherwise, since permutation $\sigma$ has no fixed points, $\beta(x) \neq \gamma(x)$ and $\lambda(\beta, \gamma, \delta)(x) = 1$ for $x\in \{0,1\}$. \sqr

\begin{figure}[!ht]
\centering
 \includegraphics[width=0.46\textwidth]{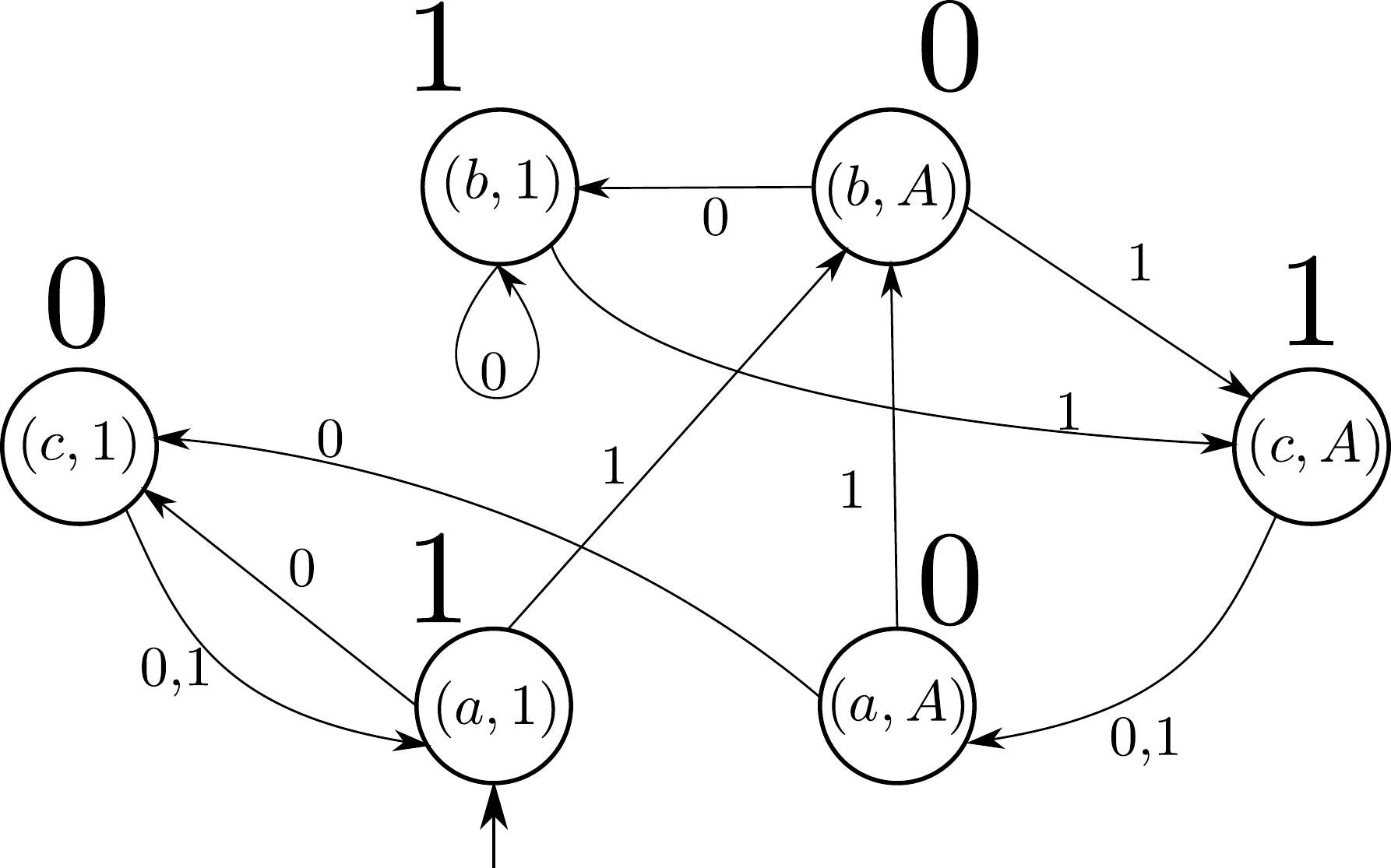}
 \caption{Automaton $L_{A,B}$ when $A$ is the adding machine and $B$ is \alesh. The output from a state is the big number next to it.}
 \label{fig:aleshd}
\end{figure}

\begin{example}\label{aleshdex}
Let $A$ be the \bb{adding machine} (see Figure \ref{fig:odom1}) with states $A$ and $\bb{1}$ (trivial state). Let \alesh{} be given by Figure \ref{fig:alesh1} and have states $\{a,b,c\}$, with initial state $a$. We consider $\Log_A F$.


Note that  
\[
A^2|_a = A|_A(a)A|_a = A,
\]
since $A|_0A|_1 = A|_1A|_0 = A$. Therefore, $A^{2^{|w|}}|_w=A$ for all $w \in X^*$ (intuitively, adding $2^n$ to a dyadic number is the same as adding $1$ to $n+1$'st digit).

We thus have $S_{A,B} \subset \{a,b,c\} \times \{A, 1\}  \times \{A\}$. Consequently, $|S_{A,B}|\leq 6$.

Let us compute the transition and the output function for $L_{A,B}$. By Proposition \ref{dmoore}, $L_{A,B}$ is a Moore machine, so  we let * stand for either $0$ or $1$ in what follows:

\begin{minipage}{0.33\textwidth}
\begin{align*}
\lambda((a,1,A),*) &= 1\\ 
\lambda((a,A,A),*) &= 0
\end{align*}
\end{minipage}
\begin{minipage}{0.33\textwidth}
\begin{align*}
\lambda((b,1,A),*) &= 1\\ 
\lambda((b,A,A),*) &= 0
\end{align*}
\end{minipage}
\begin{minipage}{0.33\textwidth}
\begin{align*}
\lambda((c,1,A),*) &= 0\\ 
\lambda((c,A,A),*) &= 1
\end{align*}
\end{minipage}
We can use this to compute the transition function:

\begin{minipage}{0.33\textwidth}
\begin{align*}
\pi((a,1,A),0) &= (c, 1, A)\\
\pi((a,1,A),1) &= (b, A, A)\\ 
\pi((a,A,A),0) &= (c, 1, A)\\ 
\pi((a,A,A),1) &= (b, A, A)
\end{align*}
\end{minipage}
\begin{minipage}{0.33\textwidth}
\begin{align*}
\pi((b,1,A),0) &= (b, 1, A)\\ 
\pi((b,1,A),1) &= (c, A, A)\\  
\pi((b,A,A),0) &= (b, 1, A)\\  
\pi((b,A,A),1) &= (c, A, A)
\end{align*}
\end{minipage}
\begin{minipage}{0.33\textwidth}
\begin{align*}
\pi((c,1,A),0) &= (a, 1, A)\\  
\pi((c,1,A),1) &= (a, 1, A)\\  
\pi((c,A,A),0) &= (a, A, A)\\   
\pi((c,A,A),1) &= (a, A, A)
\end{align*}
\end{minipage}

Since $\delta = A$ for all $(\beta, \gamma, \delta) \in S_{A,B}$, we omit it and write $(\beta, \gamma)$ for $(\beta,\gamma, A)$ 
in $\L_{A,B}$. The automaton $\L_{A,B}$ we have computed here is in in Figure \ref{fig:aleshd}.
\tri
\end{example}

Example \ref{aleshdex} calls for a more efficient notation in the case when $A$ is the adding machine and $B$ is invertible:

\begin{corollary}\label{dodomautom}
Let $A$ be the \bb{adding machine} given by automaton of Figure \ref{fig:odom1}, and assume $B$ is invertible. Then $\delta = A$ for all $(\beta, \gamma, \delta)$ in the connected component of $(B, \bb{1}, A)$ in $L_{A,B}$, and so can be omitted. After relabeling $(\beta, \gamma, \delta) \-> (\beta, \gamma)$ in $\L_{A,B}$ , we obtain the Moore machine $\hat \L_{A,B}$ with initial state $(B, \bb{1})$, 
and transition and output functions $\pi$ and $\lambda$ as specified in Table \ref{tab:logautomtablebinv}.
\end{corollary}

\bb{Note:} $\hat \L$ and $\L$ are equivalent automata: $\L(w) = \hat \L(w)$ for all words $w$.

\begin{table}[H]
\centering
\begin{tabular}{r|c|c}
& $\beta$ and $\gamma$ are both active& Exactly one of $\beta$ and $\gamma$ \\ 
& or both passive & is active\\ \hline
$\pi((\beta, \gamma), a)$ & $(\beta|_a, \gamma|_a) $ & $(\beta|_a,  (\gamma A)|_a)$\\
$\lambda((\beta,\gamma))$ & $0$ & $1$
\end{tabular}
\caption{Transition and output functions of the automaton computing $\Log_A(B)$ when $A$ is the adding machine and $B$ is invertible}
\label{tab:logautomtablebinv}
\end{table}

\pf Observe  that 
\[
A^2|_a = A|_A(a)A|_a = A,
\]
since $A|_0A|_1 = A|_1A|_0 = A$. Since the initial state is $(B,\bb{1},A)$, it follows that the rest of the states in the connected component of $L_{A,B}$ containing the initial state are of the form $(\beta, \gamma, A)$. Similarly, $\gamma \in \{\bb{1}, A\}$.

The rest follows from the construction \ref{dautom} and Prop. \ref{dmoore}. Note that $\beta(x)=\gamma(x)$ for $x \in X=\{0,1\}$ if and only if  $\beta$ and $\gamma$ are both active or both passive. \sqr

When $A$, $B$ are as in the Corollary above, it is easy to construct $\L_{A,B}$, since once can see $\beta$, $\gamma$ are active or passive by examining the diagram of the automatons $B$ and $A$. 

\begin{remark}
When $B$ is invertible, and $\beta \in S(B)$ is a state of $B$, the transition function $\lambda$ of $B$ at $\beta$, $\lambda_\beta$, takes values in $\op{Perm}(X) = \{\bb{1}, \sigma\}$.
The Table \ref{tab:logautomtablebinv} of Proposition \ref{dodomautom} can be rewritten out explicitly as Table \ref{tab:logautomtablebinv_explicit}. 
\end{remark}

\begin{table}[H]
\centering
\begin{tabular}{|c|c|c|c|c|}
\hline
$\lambda_\beta$ & $\gamma$ & $x$ & $\pi((\beta, \gamma), x)$ & 
$\lambda((\beta, \gamma), x)$ \\ \hline
$\bb{1}$ & $\bb{1}$ & $0$ & $(\pi(\beta, 0), \bb{1})$ & $0$\\
$\bb{1}$ & $\bb{1}$ & $1$ & $(\pi(\beta, 1), \bb{1})$ & $0$\\
$\sigma$ & $A$ & $0$ & $(\pi(\beta, 0), \bb{1})$ & $0$\\
$\sigma$ & $A$ & $1$ & $(\pi(\beta, 1), A)$ & $0$\\
$\bb{1}$ & $A$ & $0$ & $(\pi(\beta, 0), A)$ & $1$\\
$\bb{1}$ & $A$ & $1$ & $(\pi(\beta, 1), A)$ & $1$\\ 
$\sigma$ & $\bb{1}$ & $0$ & $(\pi(\beta, 0), \bb{1})$ & $1$\\
$\sigma$ & $\bb{1}$ & $1$ & $(\pi(\beta, 1), A)$ & $1$\\ \hline
\end{tabular}
\caption{Table \ref{tab:logautomtablebinv} with explicit values of $\pi$, $\gamma$}
\label{tab:logautomtablebinv_explicit}
\end{table}

\begin{figure}[ht]
	\begin{subfigure}{0.3\textwidth}
		\centering	
 \includegraphics[width=0.99\textwidth]{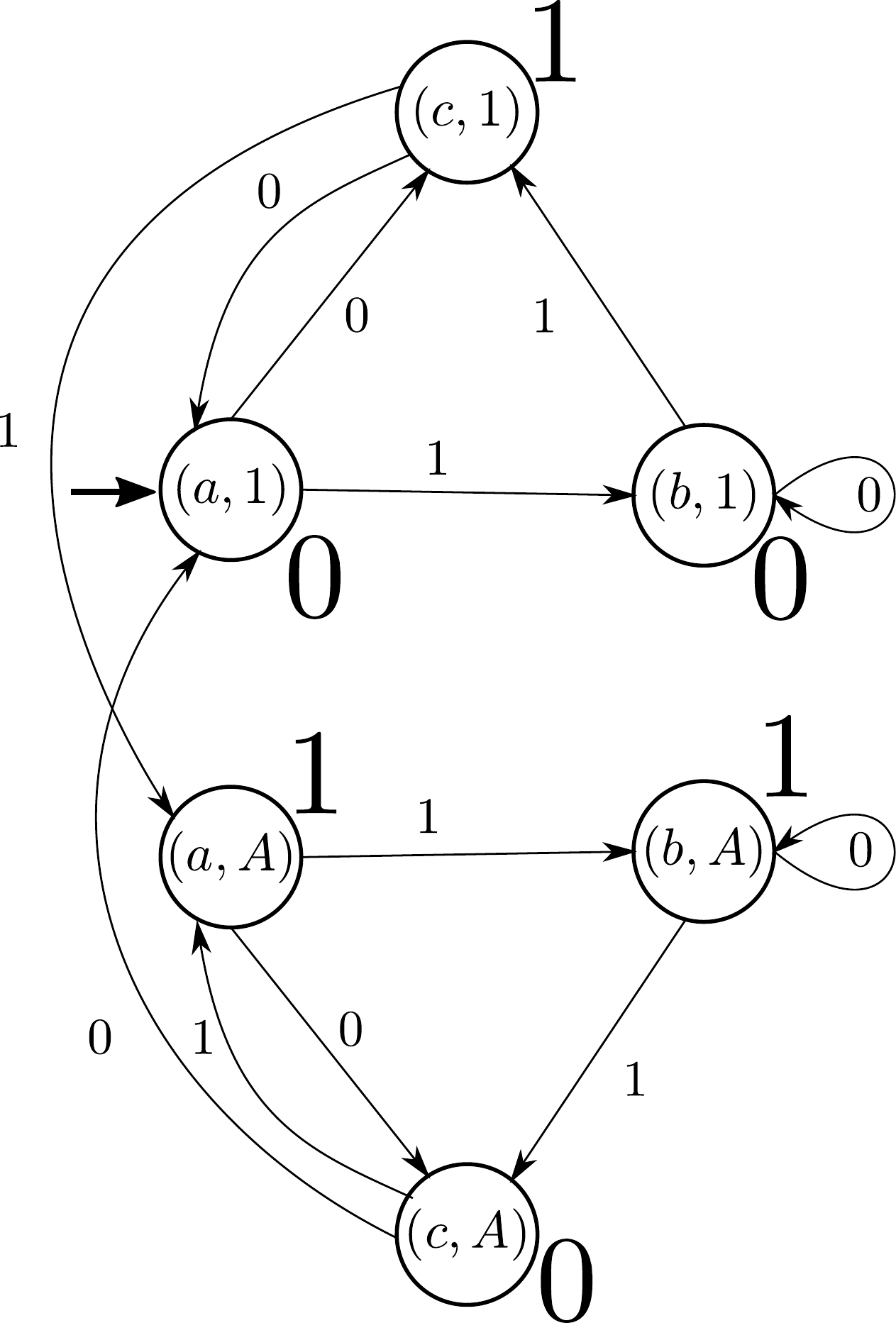}
 \caption{$L_{A,B}$}
 \label{fig:dbellaterra}
	\end{subfigure}
	\begin{subfigure}{0.6\textwidth}
		\centering	
 \includegraphics[width=0.99\textwidth]{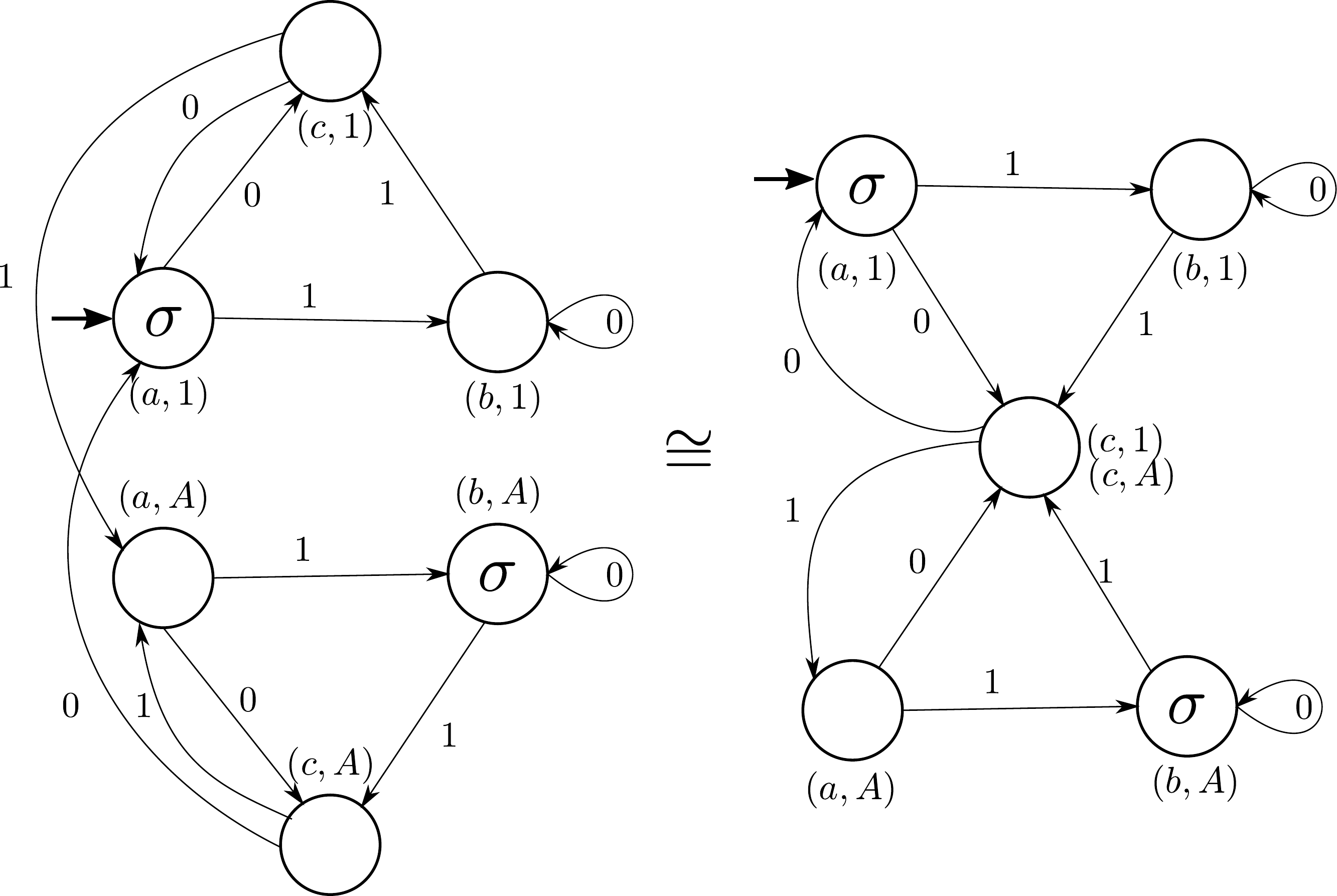}
 \caption{The delayed automaton $\sigma L_{A,B}$ is invertible}
 \label{fig:bellaterradel}
	\end{subfigure}
	\caption{Construction of $\L_{A,B}$ where $A$ is the adding machine and $B$ is Bellaterra}
\end{figure}

\begin{example}
We compute the distance automaton when $A$ is the adding machine, and $B$ is the Bellaterra automaton (Figure \ref{fig:bellaterra}). This automaton is so called because it was was studied during the summer school in Automata Groups in the Autonomous University of Barcelona in Bellaterra. An interesting property of it is that the group generated by its states is a free product of 3 copies of $\zz/2\zz$ \cite{SavVor}.

\begin{figure}[ht]
\centering
 \includegraphics[width=0.33\textwidth]{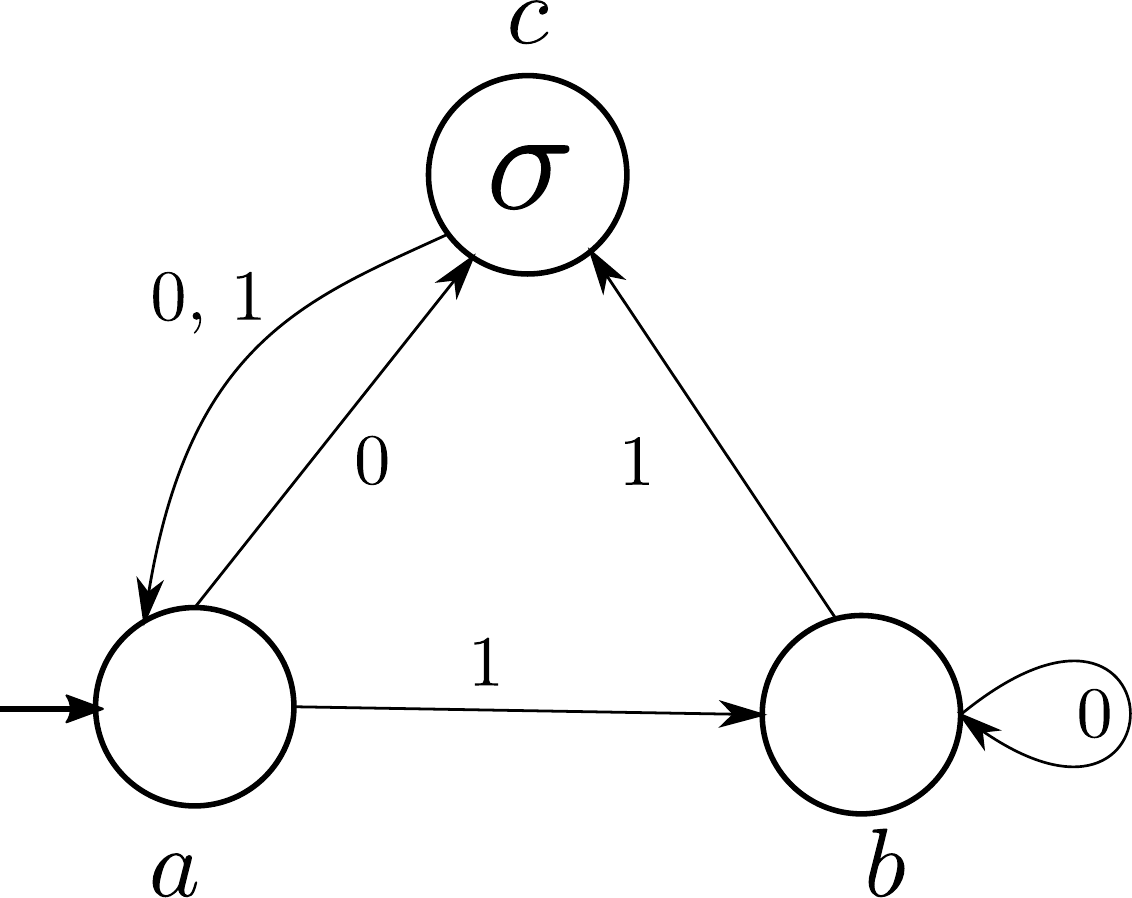}
 \caption{Bellaterra automaton}
 \label{fig:bellaterra}
\end{figure}

Using the new notation:

\begin{minipage}{0.33\textwidth}
\begin{align*}
\lambda((a,1)) &= 0\\ 
\lambda((a,A)) &= 1
\end{align*}
\end{minipage}
\begin{minipage}{0.33\textwidth}
\begin{align*}
\lambda((b,1)) &= 0\\ 
\lambda((b,A)) &= 1
\end{align*}
\end{minipage}
\begin{minipage}{0.33\textwidth}
\begin{align*}
\lambda((c,1)) &= 1\\ 
\lambda((c,A)) &= 0
\end{align*}
\end{minipage}

\begin{minipage}{0.33\textwidth}
\begin{align*}
\pi((a,1),0) &= (c, 1)\\
\pi((a,1),1) &= (b, 1)\\ 
\pi((a,A),0) &= (c, A)\\ 
\pi((a,A),1) &= (b, A)
\end{align*}
\end{minipage}
\begin{minipage}{0.33\textwidth}
\begin{align*}
\pi((b,1),0) &= (b, 1)\\ 
\pi((b,1),1) &= (c, 1)\\  
\pi((b,A),0) &= (b, A)\\  
\pi((b,A),1) &= (c, A)
\end{align*}
\end{minipage}
\begin{minipage}{0.33\textwidth}
\begin{align*}
\pi((c,1),0) &= (a, 1)\\  
\pi((c,1),1) &= (a, A)\\  
\pi((c,A),0) &= (a, 1)\\   
\pi((c,A),1) &= (a, A)
\end{align*}
\end{minipage}

In the above example, we have constructed the automaton $\tilde L$ in Figure \ref{fig:dbellaterra}. The automaton appearing in Figure \ref{fig:bellaterradel}
will be explained later.
\tri
\end{example}

\section{Automatic Exp and Logarithm of Products}
It is worthwhile to consider the operation opposite to constructing $\Log_A(B)$.

Let $|X|=2$, and let $\psi : X^\nn \-> \zz_2$ be the function that naturally identifies words in $X$ with dyadic integers:
\[
\psi(w) = \sum_{i=0}^{|w|-1} w_i 2^i.
\]

\begin{proposition}
Let $A$, $B$ be tree endomorphisms. Define a function $\Exp_A(B)_n$ on words $w$ of length $n$ by 
\[
Exp_A(B)_n(w) = A^{\psi(B(w))}(w).
\]
Then for all $n$, $\Exp_A(B)_n$ is an endomorphism of finite trees. 

The endomorphisms of finite trees $\Exp_A(B)_m$, $\Exp_A(B)_k$ agree on the levels $1,2,\dots,\min(m,k)$ on which they are both defined.
\end{proposition}
\pf We need to show that if $w \in X^*$ and $x, y \in X$, then $Exp_A(B)(wx)$ and $Exp_A(B)(wy)$ only differ in their last symbol. 
Let $n = |w|$. Observe that $\psi(B(wx))$ differs from $\psi(B(w))$ by a multiple of $2^n$. Now
\begin{align}\label{poweqn}
\Exp_A(B)_{n+1}(wx) = A^{a\cdot 2^n}A^{\psi(B(w))}(wx),
\end{align}  
where $a \in \{0,1\}$ is given by $a = B|_w(x)$. 
Note that $A^{2^n}(v) = v$ for any word $v$ of length $n$ because the length of any orbit of $A$ on the level $n$ is a factor of $2^n$.  Therefore the prefix of length $n$ of $\Exp_A(B)_{n+1}(wx)$ is given by $A^{\psi(B(w))}(w)$, i.e., it does not depend on $x$.  Thus $\Exp_A(B)$ is an endomorphism. 

The above argument also shows that $\Exp_A(B)_{n+1}$ and $\Exp_A(B)_{n}$ agree on levels $1,2,\dots,n$.  This completes the proof.
\sqr

\begin{definition}
Let $\Exp_A(B)$ denote the extension of the maps $\Exp_A(B)_n$ to the boundary of the tree $\partial \mathcal T$. We shall use the same notation for action on finite words.
\end{definition}

\begin{proposition}
Let $A$ be a tree endomorphism, and $B$ be a Moore machine. Then $\Exp_A(B)$ is an automorphism of the tree $\mathcal T$.
\end{proposition}

\pf It suffices to show that $\Exp_A(B)_n$ is invertible for all $n$.
Consider an arbitrary word $w \in X^*$ of length $n$ and a letter $x \in X$.
Recall that $B$ being a Moore machine means that $B|_w$ is constant on $X$. Therefore the value of $\Exp_A(B)_{n+1}(wx)$ (see equation (\ref{poweqn})) is given by a power of $A$ that does not depend on $x$. By assumption, $A$ is invertible, so its sections are invertible as well.  Hence $\Exp_A(B)_{n+1}|_w$ acts as a permutation of $X$. 

The proposition follows by induction on $n$. $\sqr$

%

\medbreak

\bb{Remark: }
We have constructed the Log automaton $\Log_{A}B$ for any invertible Mealy machine $B$ and any level-transitive automaton $A$ of bounded activity. By construction, the Log automaton of an invertible automaton is a Moore machine.

Therefore every invertible automaton $B$ can be written in  the form $B=\Exp_AM$, where $A$ is the adding machine (or any bounded-activity, level-transitive automaton), and $M$ is a Moore machine. Note that, in general, one cannot construct a Moore machine (synchronously) equivalent to a given Mealy machine. This construction provides an alternative.

\subsection{Logarithm of product}

\begin{proposition}\label{uniqt}
Let $A$ and $B$ be finite state automata.  Then
\[
\Exp_A(\Log_A(B)) = B
\]
as endomorphisms of the tree $\mathcal T$.  In particular, the Automatic Logarithm, as an inverse of $\Exp$, is unique. That is, if $\Exp_A B_1 = \Exp_A B_2$, then $B_1 = B_2$ as endomorphisms of trees.
\end{proposition}
\pf $\Exp_A(\Log_A(B)) = B$ by construction.
 If $\Exp_A B_1 = \Exp_A B_2$, then for any word $w$ we have $\psi(B_1(w)) = \psi(B_2(w))  \mod 2^{|w|}$. This implies $B_1(w) = B_2(w)$. \sqr

We can now argue about $\Log$ using $\Exp$. To proceed, we define:

\begin{definition}
Let $A=(S_A,\pi_A, \lambda_A, S_{A_0})$ and $B=(S_B,\pi_B, \lambda_B, S_{B_0})$ be finite automata. The \bb{sum automaton} $A \oplus B$ is the automaton with the set of states $S = S_A \times S_B \times \{0,1\}$, and transition map $\pi$ and  output map  $\lambda$ given by
\begin{align*}
\pi((s, t, c), x) &= (\pi_A(s, x), \pi_B(t, x), d), \mbox{ where}\\
d &= \begin{cases}
1 \mbox{ if } \lambda(s, x) + \lambda(t, x) + c \geq 2,\\
0 \mbox{ otherwise};
\end{cases}\\
\lambda((s, t, c), x) &= \lambda(s, x) + \lambda(t, x) + c \mod 2.
\end{align*}
\end{definition}

For a finite word $w$, the sum automaton $A \oplus B$ outputs $\psi(A(w)) + \psi(B(w))$ as a dyadic integer. The third component of a state can be understood as the carry bit.

This definition allows us to compute the Log automaton of a product.

\begin{proposition}
Let $B, C$ be invertible finite automata and $A$ be bounded-activity, level-transitive automaton. Then
\[
\Log_A(BC) = \left(\left(\Log_A B\right) C\right) \oplus \Log_A C.
\]
\end{proposition}
\pf Let $\Log_A B = a$ and $\Log_A C = c$. Then
\begin{align*}
C(w) &= A^{\psi c(w)}(w)\\
BC(w) &= A^{\psi bC(w)}(C(w))\\
      &= A^{\psi bC(w)}(A^{\psi c(w)}(w))\\
      &= A^{\psi bC(w) + \psi c(w)}(w)\\
      &= A^{\psi ((bC) \oplus c)(w)}(w)\\
      &= \Exp_A((bC) \oplus c)(w).
\end{align*}
Therefore, by Proposition \ref{uniqt},
\[
\Log_A(BC) = (bC) \oplus c,
\] which completes the proof. \sqr

\section{Distribution of lengths of chords}\label{lengthcords_thesissection}
We now approach the main goal of our investigation.
The measure we are interested in is $\mu = \mu_{A,B} := \Log_A(B)_*\nu$, where $\nu$ is the uniform Bernoulli measure on $\mathcal T$.

This measure gives the distribution of the displacement function: $d$ -- a finite integer written in binary as $w=w_0\ldots w_{n-1}$ (and thus interpreted as an element of $\zz/2^n\zz$), 
\[
\mu(\cyl{w}) = |\{v \in X^n\ :\ \Log_{A,n}(B)(v) = w\}|.
\]

We introduce this measure with the goal of studying the properties of the graphs of action, such as their diameter. For example, Pak and Malyshev prove in \cite{PakMal} that the diamter of the graph of action of the states of \alesh{} on level $n$ grows at a rate of $O(n^2)$. However computer experiments give hope that this bound can be improved to $O(n)$. Finding the connections between the measure $\mu$ and the properties of the graphs nevertheless remains an open problem.

Figure \ref{fig:schreiergraphexamples} illustrates the graphs of action with the cycle generated by the adding machine $\A$ put on a circle, and the  edges corresponding to the action of another automaton being chords in that circle, motivating the title of this section. The graph on the right has a smaller diameter.

We now proceed to examine interesting properties of $\mu$, and answer questions about it: what kind of measure is $\mu$? Is it Markov, for example?

In fact, there is an easy sufficient condition for $\mu$ to be not only Markov, but uniform Bernoulli on a cylinder. To state it, we need to make several definitions:

\begin{definition} $\sigma : X^\nn \-> X^\nn$ is the (left) \bb{shift}, defined by $\sigma(aw) = w$ for $a \in X$ and $w\in X^\nn$. We define the left shift $\sigma : X^* \-> X^*$ for finite words $w$ in the same way.
\end{definition}

\begin{definition}\label{delautom}
When $\L$ is a Moore machine, the \bb{delayed automaton} $\sigma \L$ is the automaton that computes the composition $\sigma \circ \L$. It has the same states, initial state and the transition function as $\L$, but the output function $\sigma\lambda$ is given by
\[
\sigma\lambda(s,x) = \lambda(\pi(s,x)),
\]
which is well-defined when $\L$ is a Moore machine. 
\end{definition}

When $\L$ is Moore, for any finite word $w \in X^*$ and $x \in X$, 
\[
\L(wx)=\L(0)\sigma \L(w)=\L(1)\sigma \L(w).
\]

\begin{proposition}\label{bernmeas}
Let $X$ be a finite alphabet. Let $\L$ be a Moore machine with initial state $s_0$, and let $a = \lambda(s_0)$. Let $\nu$ be the uniform Bernoulli measure on $X^\nn$.

Then $\mu=\L_*\nu$ is supported on the cylinder $aX^\nn$, and $\mu|_a = (\sigma \L)_*\nu$. If $\sigma L$ is invertible, $\mu|_a$ is uniform Bernoulli (i.e. $\mu|_a=\nu$).
\end{proposition}
\pf First, note that 
\[
\mu(\cyl a)  = \L_*\nu(\cyl a) = \nu(\L^{-1}(\cyl a)) = \nu(\lambda_{s_0}^{-1}(a)\cyl{}) = \nu(\cyl{}) = 1.
\]
Now $\mu_a=(\L_*\nu)|_a=(\sigma L)_*\nu$, since for all $v \in X^*$,
\begin{align*}
(\sigma L)_*\nu(\cyl{v}) &= \nu\lp (\sigma L)^{-1}(\cyl{v})\rp\\
&=  \nu\lp  \L^{-1} (\sigma^{-1}\lp(\cyl{v})\rp \rp\\
&=  \nu\lp  \L^{-1} \lp \bigsqcup_{x \in X} x\cyl{v}\rp \rp\\
&=  \nu\lp  \L^{-1} (a\cyl{v})\rp\\
&= \L_*v\lp a\cyl{v}\rp\\
&=(\L_*v)|_a(\cyl{v}) &\mbox{(since, as noted, $\L_*\nu(\cyl{a})=1$)}.
\end{align*}
Thus $(\sigma \L)_*\nu = (\L_*\nu)|_a$.

If $\sigma \L$ is invertible, then $(\sigma \L)_*\nu = \nu$ by Proposition \ref{invpush}. This completes the proof. \sqr

\begin{corollary}Let $X$, $A$, $B$ and $\L=\L_{A,B}$ be as in Prop \ref{dautom} (so $B$ is invertible, and $\L=\L_{A,B}$ is Moore). Let $\nu$ be the uniform Bernoulli measure on $X^*$.

Then $\mu=\Log_A(B)_*\nu$ is supported on $\cyl{\L(0)}$, and $\mu|_{\L(0}=\nu$.
\end{corollary}

\begin{example}\label{bellaex}
Let $\L=\L_{A,B}$ with $A$ being the adding machine, and $B$ being the Bellaterra automaton defined in Fig. \ref{fig:bellaterra}. Then $\L(0)=\L(1)=0$. The delayed automaton $\sigma \L$ is shown in Figure \ref{fig:bellaterradel}, and it is invertible (but not minimal: can be reduced to an automaton with $5$ states).

Therefore, $\mu=\Log_A(B)_*\nu$ is the uniform Bernoulli measure supported on $\cyl{0}$, i.e. $\mu|_0=\nu$ and $\mu|_1 = 0$.
\tri
\end{example}

\begin{figure}[ht]
\begin{subfigure}{0.6\textwidth}
 \centering
 \includegraphics[width=0.9\textwidth]{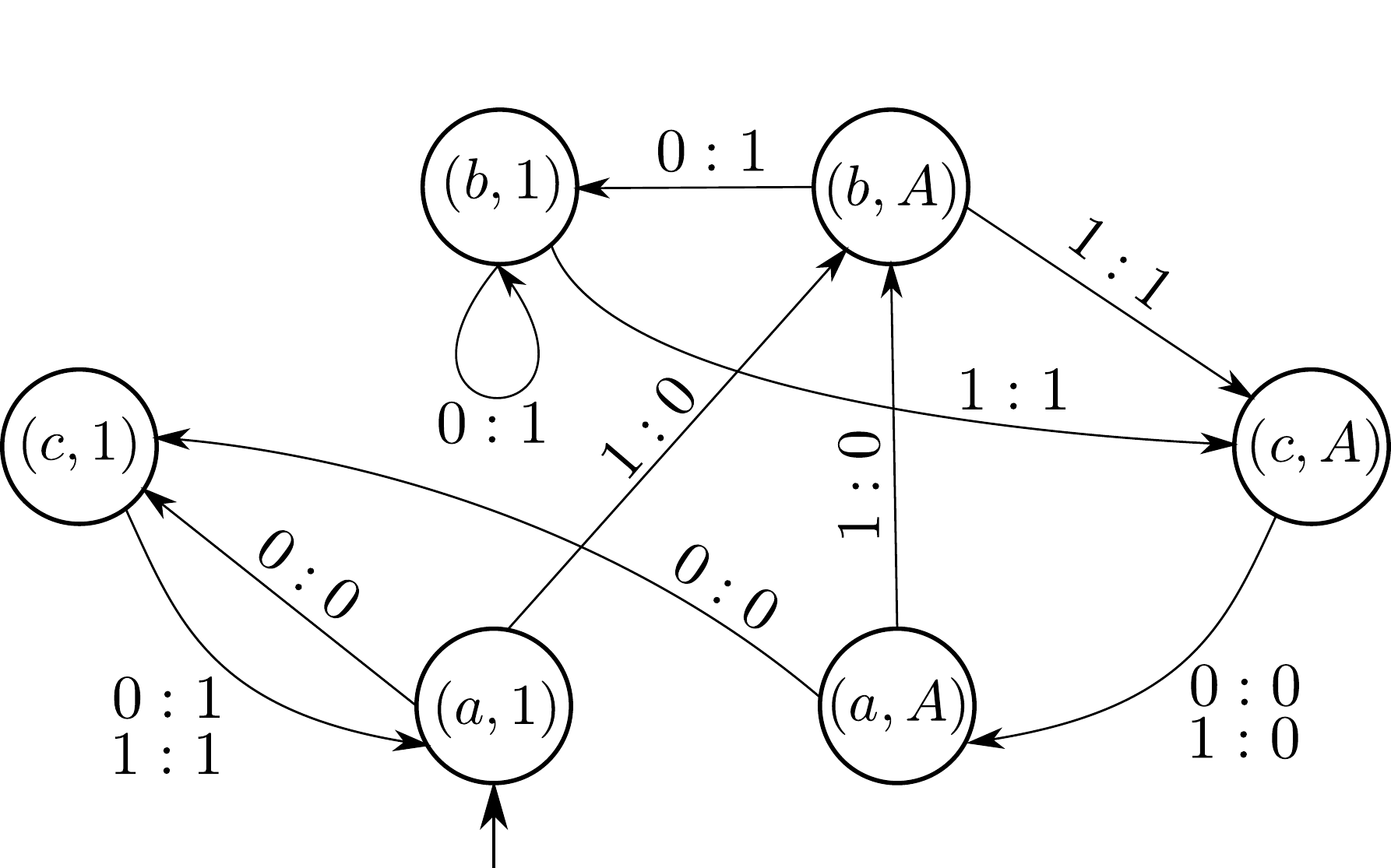}
 \caption{$\sigma \L_{A,B}$}
 \label{fig:aleshddel}
\end{subfigure}
\begin{subfigure}{0.4\textwidth}
 \centering
 \includegraphics[width=0.9\textwidth]{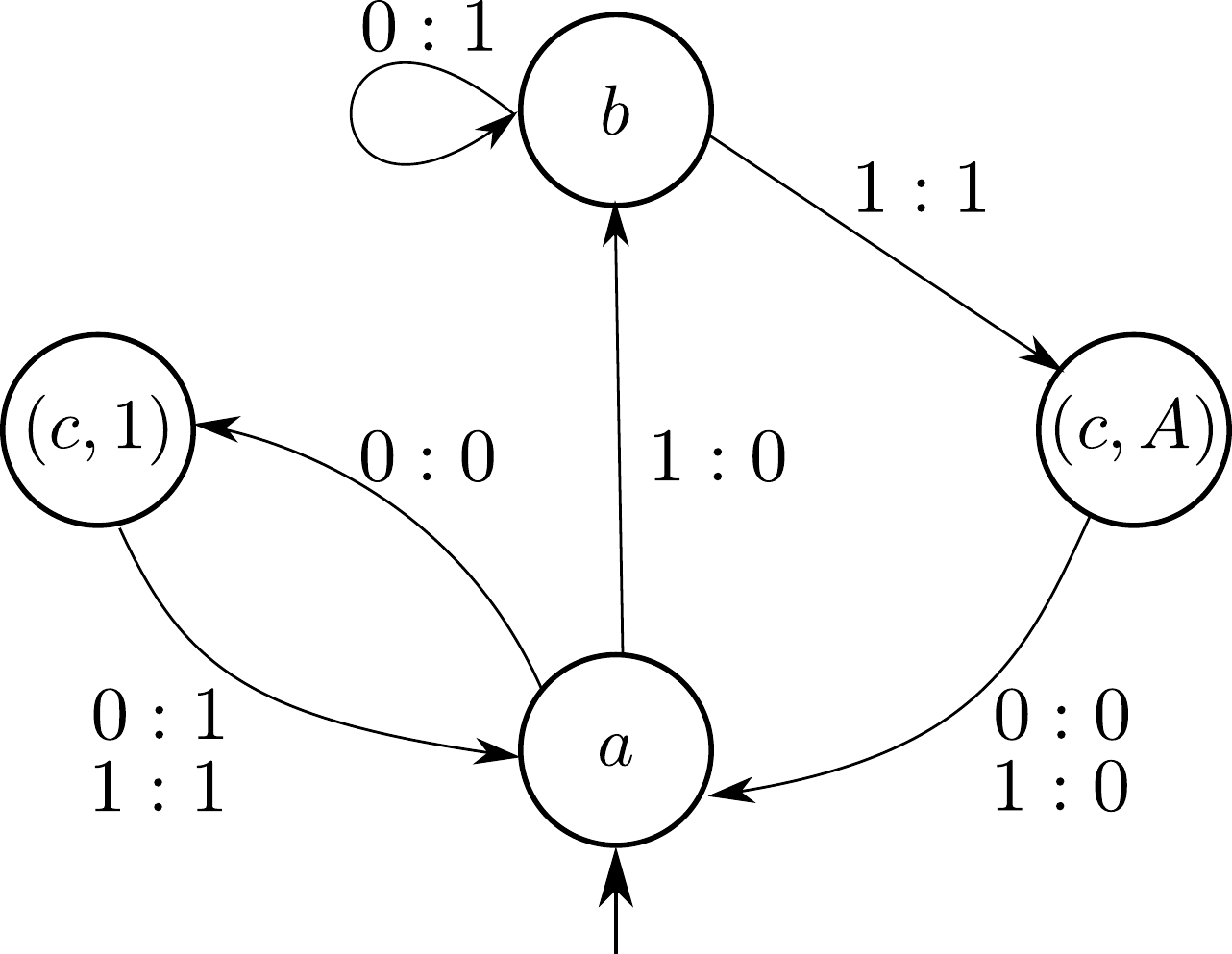}
 \caption{$\sigma \L_{A,B}$ minimized.}
 \label{fig:aleshdmin}
\end{subfigure}
\caption{Automatons $\sigma \L_{A,B}$ and its minimization}
\end{figure}

\begin{figure}[ht]
\centering
\includegraphics[width=0.6\textwidth]{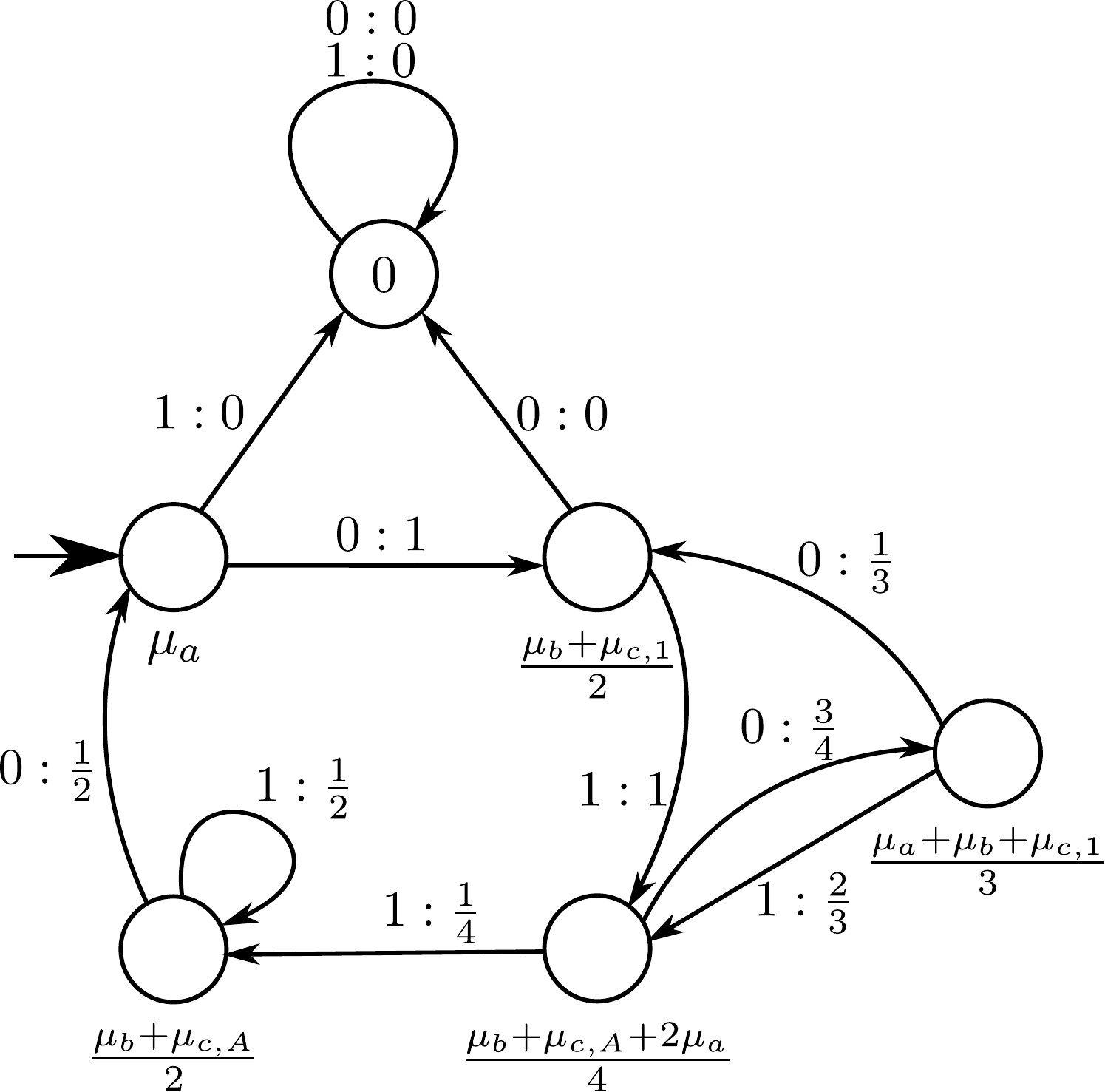}
\caption{Automaton that computes $\mu_{A, B}|_1$ for $A$ the adding machine and $B$ -- \alesh, defined in Figure \ref{fig:alesh1}}
\label{fig:aleshmeas}
\end{figure}

Proposition \ref{bernmeas} demonstrates that when $B$ is invertible, the delayed automaton $\sigma \L_{A,B}$ can be useful for examining $\mu_{A,B}$. We make use of it again for what follows:

\begin{theorem}\label{aleshmex}
Let $\nu$ be the uniform Bernoulli measure. In the case $A$ is the adding machine and $B$ is \alesh{} (see Figure \ref{fig:alesh1}), the measure $\mu_{A,B}=\Log_A(B)_* \nu$ is finite-state. Furthermore, $\mu_{A,B}|_0 = 0$, and automaton in Figure \ref{fig:aleshmeas} computes $\mu_{A,B}|_1$ (in the sense of Definition \ref{fsmeasdefin}).
\end{theorem}
\pf Write $\mu=\mu_{A,B}$. By Proposition \ref{bernmeas} and the already computed $\L=\L_{A,B}$ in Fig \ref{fig:aleshd}, $\mu|_0 = 0$, and 
the measure is supported on the cylinder $\cyl{1}$, with $\mu|_1 = (\sigma \L)_*\nu$. We thus point our attention to $\sigma \L$, shown  in Figure \ref{fig:aleshddel}.

First, observe that the automaton $\sigma \L$ is not minimal. After identifying states $(a,1)$ and $(a,A)$ into state $a$, and identifying states $(b, 1)$ and $(b, A)$ into state $b$, we obtain a minimal automaton $\L$ (Figure \ref{fig:aleshdmin}).

Recall that the states $a$, $b$, $(c,1)$ and $(c,A)$ of the automaton $\sigma\L_{A,B}$ are sections of $\L_{A,B}$, and can be seen as endomorphisms whose automata coincide with $\L_{A,B}$ except for the initial state (see Remark \ref{statesectremark}), i.e., $\L_{A,B} = a$, $\L_{A,B}|_1 = b$, etc. 

If $g$ is an action on the tree $\mathcal T$, we write $\mu_g$ for $g_*\nu$. Thus we are interested in $\mu_a = \mu = a_*\nu = (\sigma \L)_*\nu$, and we compute it by writing down its sections in terms of $\mu_a$, $\mu_b$, $\mu_{c,1}$ and $\mu_{c,A}$.

We apply Corollary \ref{imsecber} to $\L_{A,B}$ to obtain the sections of $\mu$ by  $x\in X=\{0,1\}$. On the right, we evaluate these measures on the cylindrical sets of the form $xX^\nn$, so that we could continue the computation by applying Proposition \ref{sumsec}.

\begin{minipage}{0.49\textwidth}
\begin{align*}
\mu_a|_0 &= \frac{\mu_b + \mu_{c,1}}{2}\\
\mu_a|_1 &= 0\\
\mu_b|_0 &= 0\\
\mu_b|_1 &= \frac{\mu_b + \mu_{c,A}}{2}\\
\mu_{c,1}|_0 &= 0\\
\mu_{c,1}|_1 &= \mu_a\\
\mu_{c,A}|_0 &= \mu_a\\
\mu_{c,A}|_1 &= 0
\end{align*} 
\end{minipage}
\begin{minipage}{0.49\textwidth}
\begin{align*}
\mu_a(\cyl 0) &= 1\\
\mu_a(\cyl 1) &= 0\\ 
\mu_b(\cyl 0) &= 0\\
\mu_b(\cyl 1) &= 1\\
\mu_{c,1}(\cyl 0) &= 0\\
\mu_{c,1}(\cyl 1) &= 1 \\
\mu_{c,A}(\cyl 0) &= 1\\
\mu_{c,A}(\cyl 1) &= 0 
\end{align*}
\end{minipage}

Having expressed the sections by one character in terms of each other, we have obtained a set of recursive relations which allows us to compute sections by arbitrary words. To find the set of all sections, we proceed by repeatedly computing sections using Proposition \ref{sumsec}. We find:

\begin{minipage}{0.49\textwidth}
\begin{align*}
\frac{\mu_b + \mu_{c,1}}{2}|_0 &= 0\\
\frac{\mu_b + \mu_{c,1}}{2}|_1 &= \frac{\mu_b + \mu_{c,A} + 2\mu_a}{4}\\
\frac{\mu_b + \mu_{c,A}}{2}|_0 &= \mu_a\\
\frac{\mu_b + \mu_{c,A}}{2}|_1 &= \frac{\mu_b + \mu_{c,A}}{2}
\end{align*} 
\end{minipage}
\begin{minipage}{0.49\textwidth}
\begin{align*}
\frac{\mu_b + \mu_{c,1}}{2}(\cyl 0) &= 0\\
\frac{\mu_b + \mu_{c,1}}{2}(\cyl 1)&= 1\\
\frac{\mu_b + \mu_{c,A}}{2}(\cyl 0) &= \frac{1}{2}\\
\frac{\mu_b + \mu_{c,A}}{2}(\cyl 1) &= \frac{1}{2}
\end{align*} 
\end{minipage}

And again:

\begin{minipage}{0.49\textwidth}
\begin{align*}
\frac{\mu_b + \mu_{c,A} + 2\mu_a}{4}|_0 &= \frac{\mu_a + \mu_b + \mu_{c,1}}{3}\\
\frac{\mu_b + \mu_{c,A} + 2\mu_a}{4}|_1 &= \frac{\mu_b + \mu_{c,A}}{2}
\end{align*} 
\end{minipage}
\begin{minipage}{0.49\textwidth}
\begin{align*}
\frac{\mu_b + \mu_{c,A} + 2\mu_a}{4}(\cyl 0) &= \frac{3}{4}\\
\frac{\mu_b + \mu_{c,A} + 2\mu_a}{4}(\cyl 1) &= \frac{1}{4}
\end{align*} 
\end{minipage}

Finally:

\begin{minipage}{0.49\textwidth}
\begin{align*}
\frac{\mu_a + \mu_b + \mu_{c,1}}{3}|_0 &= \frac{\mu_b + \mu_{c,1}}{2}\\
\frac{\mu_a + \mu_b + \mu_{c,1}}{3}|_1 &= \frac{\mu_b + \mu_{c,A} + 2 \mu_a}{4}
\end{align*} 
\end{minipage}
\begin{minipage}{0.49\textwidth}
\begin{align*}
\frac{\mu_a + \mu_b + \mu_{c,1}}{3}(\cyl 0) &= \frac{2}{3}\\
\frac{\mu_a + \mu_b + \mu_{c,1}}{3}(\cyl 1) &= \frac{1}{3}
\end{align*} 
\end{minipage}

Since we have obtained no new sections at this step, the sections so far are all the sections of $\mu$. We have all the data now to build the automaton in Figure \ref{fig:aleshmeas} that computes $\mu|_1$. \sqr

The preceding example shows that $\mu_{A,B}$ is finite-state (in the sense of Definition \ref{fsmeasdefin}) in the case when $A$ is the adding machine and $B$ is \alesh.
It should be noted that for some choices of automaton $B$ the measure $\mu_{A,B}$ is not finite-state.

\begin{example}\label{lampex} Let $A$ be the adding machine and $B$ be the Lamplighter automaton; see Figure \ref{fig:lamplighter}.
Then the measure $\mu_{A,B}$ is not finite-state as shown below. \tri
\end{example}

\begin{figure}[ht]
\centering
\includegraphics[width=0.35\textwidth]{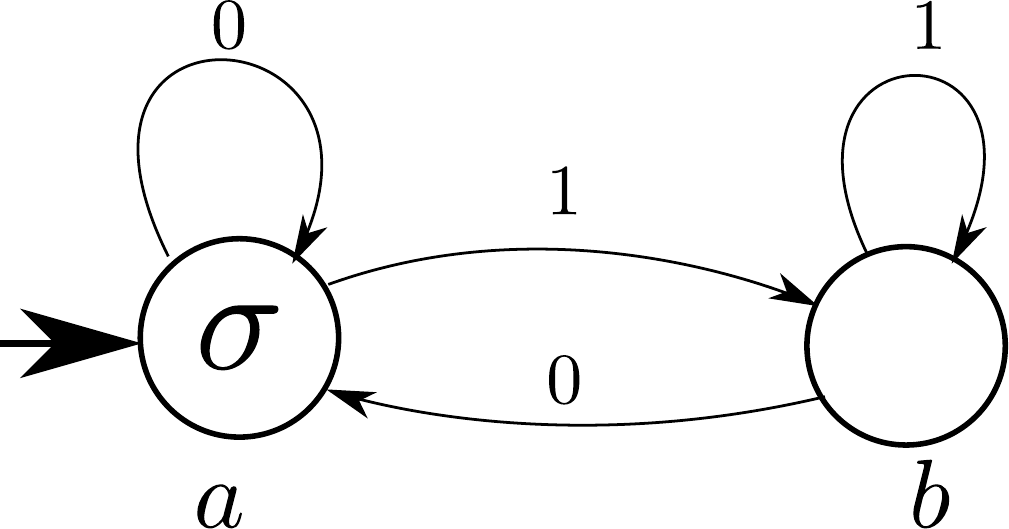}
\caption{The lamplighter automaton}
\label{fig:lamplighter}
\end{figure}

We compute the automaton $\L_{A,B}$ using Theorem \ref{dautom}:

\begin{minipage}{0.2\textwidth}
\begin{align*}
\lambda((a,1)) &= 1\\ 
\lambda((a,A)) &= 0
\end{align*}
\end{minipage}
\begin{minipage}{0.2\textwidth}
\begin{align*}
\lambda((b,1)) &= 0\\ 
\lambda((b,A)) &= 1
\end{align*}
\end{minipage}
\begin{minipage}{0.3\textwidth}
\begin{align*}
\pi((a,1),0) &= (a, 1)\\
\pi((a,1),1) &= (b, A)\\ 
\pi((a,A),0) &= (a, 1)\\ 
\pi((a,A),1) &= (b, A)
\end{align*}
\end{minipage}
\begin{minipage}{0.3\textwidth}
\begin{align*}
\pi((b,1),0) &= (a, 1)\\ 
\pi((b,1),1) &= (b, 1)\\  
\pi((b,A),0) &= (a, A)\\  
\pi((b,A),1) &= (b, A)
\end{align*}
\end{minipage}
\vspace{0.1in}

\begin{figure}[ht]
\centering
\begin{subfigure}{0.35 \textwidth}
\centering
\includegraphics[width=0.99\textwidth]{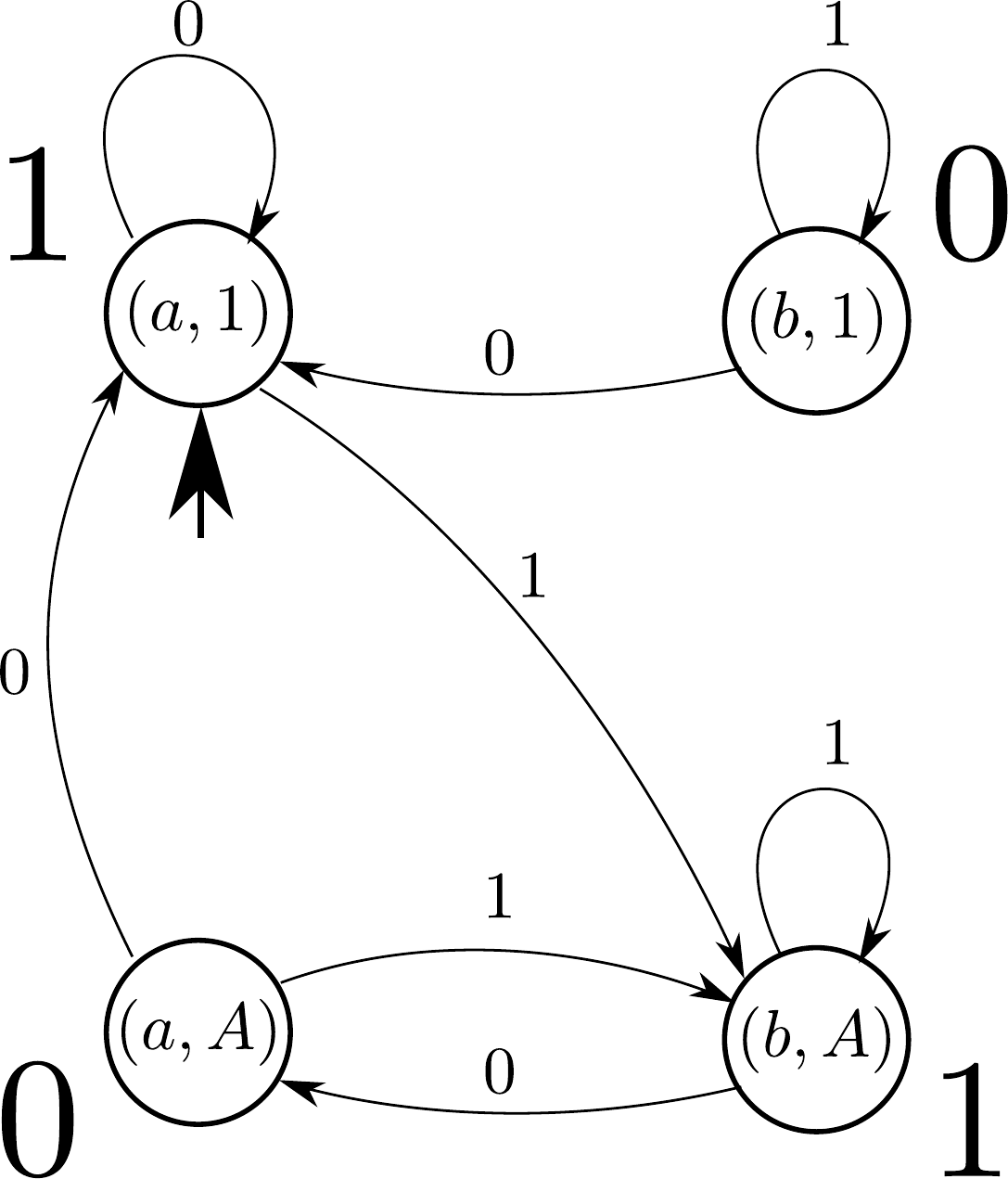}
\caption{$\L_{A,B}$}.
\label{fig:dlamplighter}
\end{subfigure}
\hspace{0.09\textwidth}
\begin{subfigure}{0.4 \textwidth}
\centering
\includegraphics[width=0.99\textwidth]{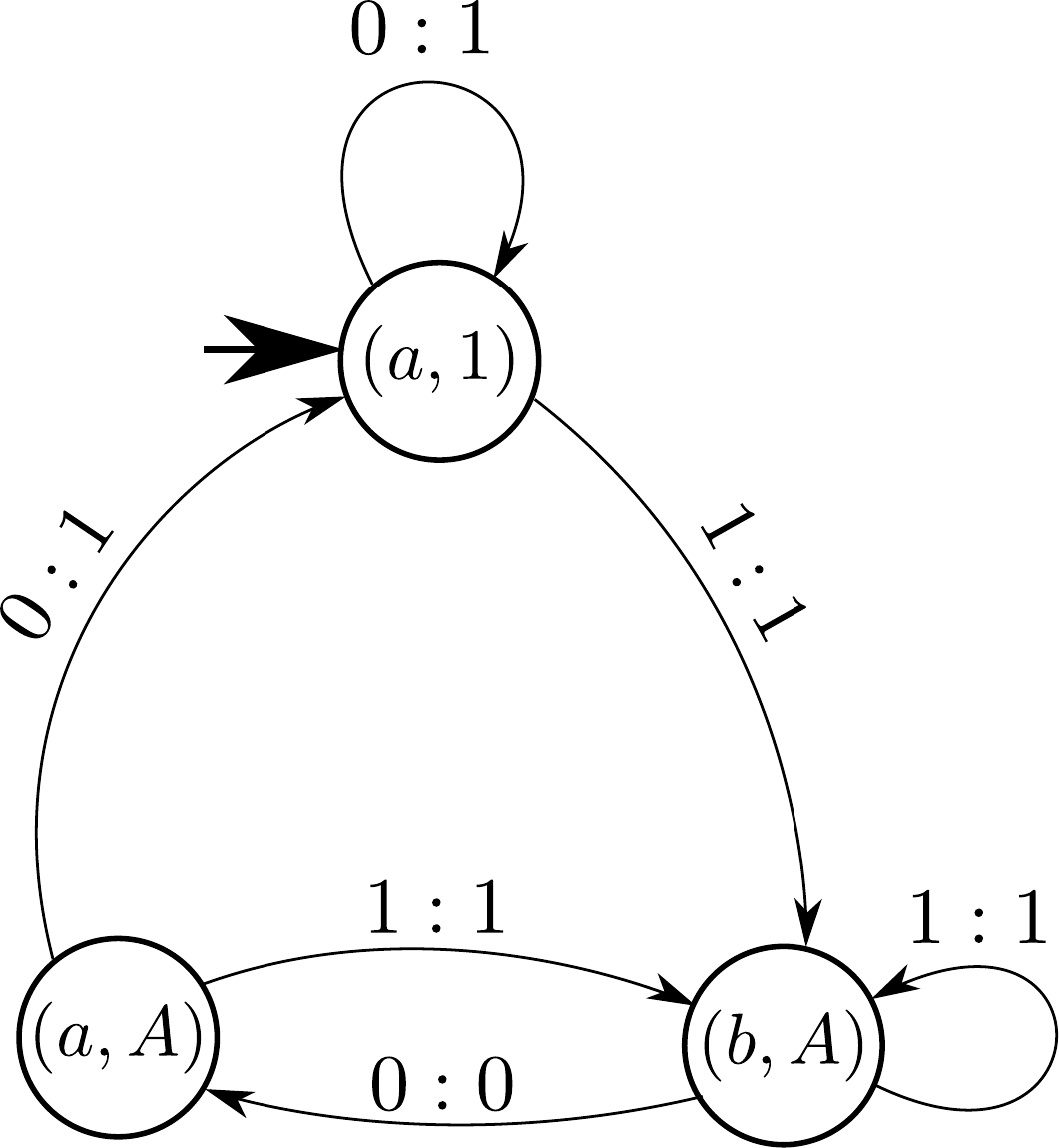}
\caption{$\sigma \L$; states unreachable from the initial state $(a,1)$ not shown}.
\label{fig:sigmadlamplighter}
\end{subfigure}

\begin{subfigure}{0.5 \textwidth}
\centering
\includegraphics[width=0.99\textwidth]{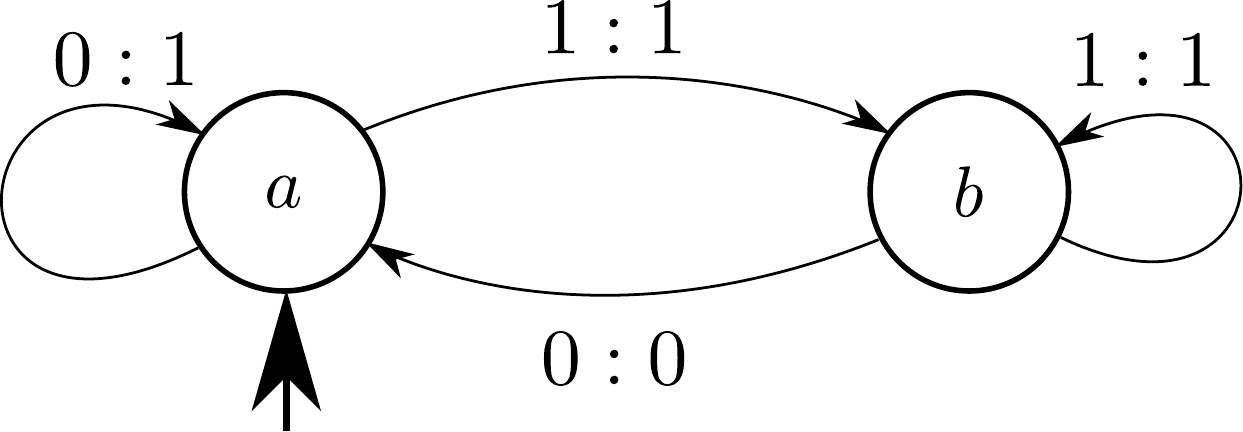}
\caption{$\sigma \L$ minimized}.
\label{fig:sigmadlamplightermin}
\end{subfigure}
\caption{$\L=\L_{A,B}$ and $\sigma \L$ for $A$ the adding machine, $B$ the Lamplighter.}
\end{figure}

The diagrams of the automata $\L_{A,B}$ and $\sigma \L_{A,B}$ are shown in Figures \ref{fig:dlamplighter} and \ref{fig:sigmadlamplighter}, respectively. 
Since $(b,1)$ is not reachable from the initial state $(a,1)$, it is omitted in Figure \ref{fig:sigmadlamplighter}. The automaton in that figure is not minimal; states $(a,1)$ and $(a,A)$ can be identified. The minimized automaton is shown in Figure \ref{fig:sigmadlamplightermin}; the relabeling is $a=(a,1)=(a,A)$, $b=(b,A)$, and $(b,1)$ is discarded as unreachable from the initial state $a$.

Noting that $\mu_{A,B}$ is supported on $\cyl{1}$ (by Proposition \ref{bernmeas}), we now point our attention to the measure $\tilde\mu=\mu_{A,B}|_1$. Using Corollary \ref{imsecber} for the minimized $\sigma \L$ in Figure \ref{fig:sigmadlamplightermin} and the notation of Example \ref{aleshmex}, we get:

\begin{minipage}{0.49\textwidth}
\begin{align*}
\mu_a|_0 &= 0\\
\mu_a|_1 &= \frac{1}{2}(\mu_a + \mu_b)\\
\mu_b|_0 &= \mu_a\\
\mu_b|_1 &= \mu_b
\end{align*} 
\end{minipage}
\begin{minipage}{0.49\textwidth}
\begin{align*}
\mu_a(\cyl 0) &= 0\\
\mu_a(\cyl 1) &= 1\\ 
\mu_b(\cyl 0) &= \frac{1}{2}\\
\mu_b(\cyl 1) &= \frac{1}{2} 
\end{align*}
\end{minipage}

Now let $\mu_0 := \mu_a$ and $\mu_n := \mu_{n-1}|_1$. Again we use Corollary \ref{sumsec}:
\begin{align*}
\mu_1 &= \frac{(\mu_a + \mu_b)}{2}\\
\mu_2 &= \mu_1|_1 =\frac{1}{2}\lp\mu_a|_1 + \frac{\mu_b|_1}{2}\rp  / \mu_1 \pcyl 1\\
&=\frac{(\mu_a + 2 \mu_b)}{4} \cdot \frac{4}{3}\\
&=\frac{(\mu_a + 2 \mu_b)}{3}\\
\mu_3 &= \mu_2|_1 =\frac{(\mu_a + 3\mu_b)}{4}\\
&\ldots\\
\mu_n &= \mu_{n-1}|_1 = \ldots
\end{align*}
All this leads to the following:
\begin{proposition} Let $\mu_0 = \mu_a$, and $\mu_n = \mu_{n-1}|_1$, for $n \in \nn$. Then
\begin{align*}
\mu_n &= \frac{\mu_a + n \mu_b}{n+1}\\
\mu_n(\cyl 0) &= \frac{ n }{2(n+1)}\\
\mu_n(\cyl 1) &= \frac{ n + 2 }{2(n+1)}
\end{align*}
\end{proposition}
\pf By induction on $n$. The proposition holds for $n=0$. Assuming it holds for $n=k$, by Corollary \ref{sumsec}:
\begin{align*}
\mu_{k+1} = \mu_k|_1 &= \frac{1}{k+1} \lp \frac{\mu_a + \mu_b}{2} + \frac{1}{2}k\mu_b\rp
 / \lp \frac{k + 2}{2(k + 1)} \rp\\
&= \frac{\mu_a + (k+1)\mu_b}{k+2}.\sqr
\end{align*}

Note that measures $\mu_n$ are all distinct. 

\begin{corollary}\label{not_fs}
$\mu_{A,B}$ is not finite-state when $A$ is the adding machine and $B$ is Lamplighter.
\end{corollary}

\begin{corollary}
$\mu_n$ for $n=0,1,2,\ldots$ are all the nontrivial sections of $\tilde\mu$. 
\end{corollary}
\pf This immediately follows from observing that $\mu_n|_0 = \mu_0$ for $n>0$:
\begin{align*}
\mu_n|_0 &= \frac{\mu_a + n\mu_b}{n+1}|_0  = \frac{1}{n+1}\frac{n \mu_a}{2}\frac{2(n+1)}{n} = \mu_a = \mu_0.
\end{align*}

\begin{figure}
\centering
\includegraphics[width=0.99\textwidth]{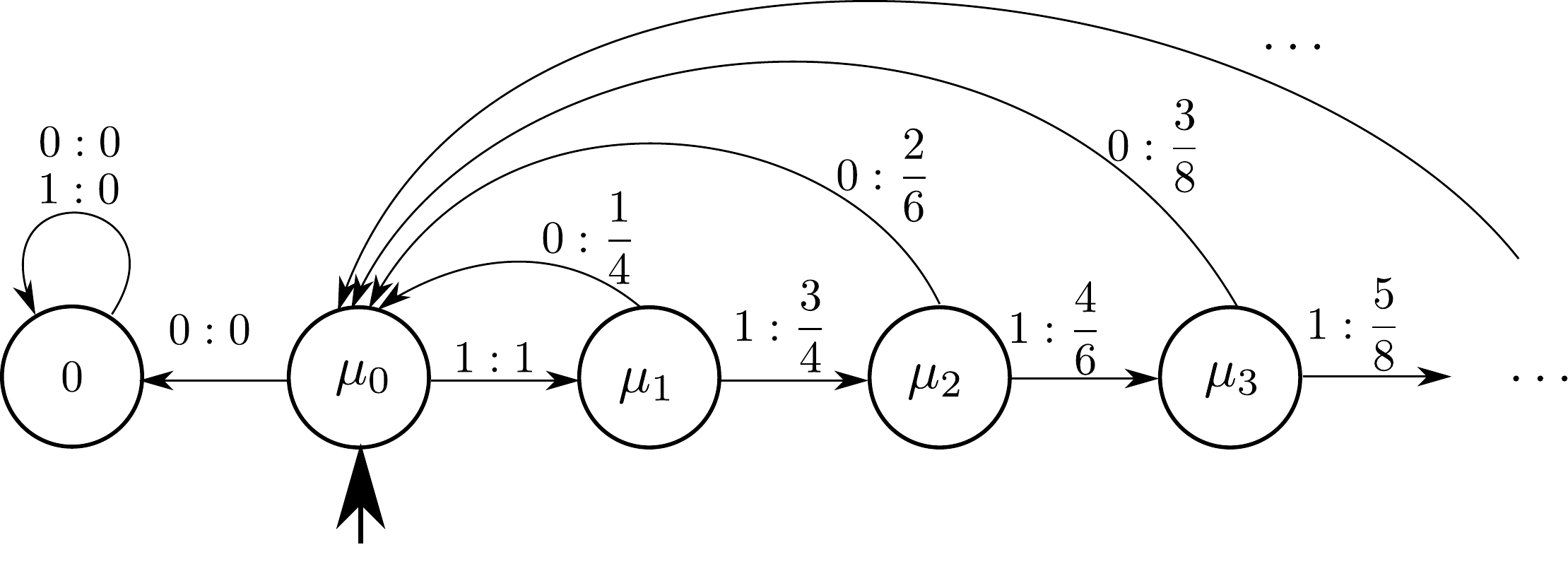}
\caption{The infinite automaton computing $\tilde\mu_{A,B}$ where $A$ is the adding machine, and $B$ is the Lamplighter automaton.}
\label{fig:lamplightermeas}
\end{figure}

The (infinite) automaton that computes $\tilde\mu$ is shown in Figure \ref{fig:lamplightermeas}.\sqr

\medskip

Observe that the computations in these examples are almost linear.  The following proposition makes this notion precise.

\begin{proposition}\label{meassec}
Let $X=\{x_0, \ldots, x_{k-1}\}$ be a finite alphabet, $L$ be a Mealy machine with states $S=\{g_0, \ldots, g_{n-1}\}$, and $\nu$ be a Bernoulli measure given by a vector $p=(p(x_0),\ldots, p(x_{k-1}))$. 
For any vector $v=(a_0, a_1, \ldots, a_{n-1})\in \rr$ let
\[
\mu_v = \sum_{i=0}^{n-1} a_i {g_i}_*\nu.
\]
Then for any $x \in X$ there exists an $n\times n$ matrix $M_x$ and an $n$-dimensional vector $p_x$ such that 
\[
\mu_v|_x = \mu_w
\]
with
\[
w = \frac{M_x v}{p_x . v}.
\]
The entries of $M_x$ and coordinates of $p_x$ are given by
\begin{align*}
M_x(i,j) &= \sum_{y: \pi(g_i, y) = g_j \mbox{ and } \lambda(g_i, y) = x} p(y);\\
p_x(j) &= \sum_{i=0}^{n-1} M_x(i,j).
\end{align*}
\end{proposition}
\pf From Proposition \ref{imsecber} and Corollary \ref{sumsec}:
\begin{align*}
\left .\lp \sum_{i=0}^n a_i {g_i}_* \nu \rp \right|_x &= \frac{\displaystyle\sum_{i=0}^n a_i {g_i}_* \nu(\cyl x) ({g_i}_*\nu)|_x}{\displaystyle\sum_{i=0}^n a_i {g_i}_* \nu(\cyl x)}\\
&=\frac{\displaystyle\sum_{i=0}^n a_i \displaystyle \sum_{y \in \lambda_{g_i}^{-1}(x)} p(y) \pi(g_i, y)_*\nu}
{\displaystyle\sum_{i=0}^n a_i \displaystyle \sum_{y \in \lambda_{g_i}^{-1}(x)} p(y)}.
\end{align*}
The proposition follows. \sqr
\begin{corollary} Let 
\[
\phi_x(v) := \frac{M_x v}{p_x \cdot v}.
\]
Then $\mu_v$ is finite-state if and only if the orbit of $v$ under the action of the semigroup generated by $\phi_x$, $x \in X$ is finite. The graph of the action is the transition diagram of the automaton that computes $\mu_{[v]}$.
\end{corollary}

The above corollary can be made simpler once we consider $v$ as an element of $\rr \pr^n$. For $v=(a_0,a_1,\ldots,a_{n-1})$, write $[v]=[a_0:a_1:\ldots:a_{n-1}]\in \rr \pr^n$ and let
\[
\mu_{[v]} := \frac{\mu_v}{\mu_v \pcyl{}}.
\]
This is well defined and
\[
[\phi_x(v)] = [M_x v].
\]
\begin{corollary} $\mu_{[v]}$ is finite-state if and only if the orbit of $[v]$ under the action of the semigroup generated by $\langle M_x : x \in X \rangle$ is finite.
\end{corollary}

In the special case when $\nu$ is the uniform Bernoulli measure, it is convenient to use matrices $\tilde M_x$ with entries
\[
\tilde M_x(i, j) = \sum_{y: \pi(g_i, y) = g_j \mbox{ and } \lambda(g_i, y) = x} 1.
\]
Similarly, set $\tilde p_x = |X|p_x$. By definition, $\tilde M_x = |X| M_x$, $[\tilde M_x v] = [M_x v]$, and $\phi_x(v)=\tilde M_x v/\tilde p_x \cdot v$.  However $\tilde M_x$ has integer entries: $M_x(i,j) \in \{0,1,\ldots, |X|\}$.
\begin{corollary}
In the case $\nu$ is uniform Bernoulli, the measure $\mu_{[v]}$ is finite-state if and only if the orbit of $[v]$ under the action of the multiplicative semigroup generated by integer matrices $\tilde M_x$, $x \in X$ is finite. 
\end{corollary}

\begin{example}
When $L=L_{A, B}$, where $A$ is the adding machine and $B$ is \alesh{} given by Figure \ref{fig:alesh1}, we have

\begin{align*}
\tilde M_0 &=\left(\!
\begin{array}{cccc}
 0 & 0 & 0 & 2 \\
 1 & 0 & 0 & 0 \\
 1 & 0 & 0 & 0 \\
 0 & 0 & 0 & 0 \\
\end{array}
\!\right)\\
\tilde p_0 &= (2,0,0,2)\\
\phi_0(v) &= {\tilde M_0 v}/{\tilde p_0 \cdot v}\\
\tilde M_1 &= \left(\!
\begin{array}{cccc}
 0 & 0 & 2 & 0 \\
 0 & 1 & 0 & 0 \\
 0 & 0 & 0 & 0 \\
 0 & 1 & 0 & 0 \\
\end{array}
\!\right)\\
p_1 &= (0,2,2,0)\\
\phi_1(v) &= {M_1 v}/{p_1\cdot v}
\end{align*}
The orbit of $(1,0,0,0)$ under the action of $\langle \phi_0, \phi_1 \rangle$ is 
\[
((0,0,0,0),(0,1/2,0,1/2),(0,1/2,1/2,0),(1/3,1/3,1/3,0),(1/2,1/4,0,1/4),(1,0,0,0)).
\]
These correspond to the states in Figure \ref{fig:aleshmeas}. 

Equivalently, the orbit of $[1:0:0:0]$ under the action of $\langle \tilde M_0, \tilde M_1 \rangle $
is 
\[
([0:0:0:0],[0:1:0:1],[0:1:1:0],[1:1:1:0],[2:1:0:1],[1:0:0:0]).
\]
\tri
\end{example}

\begin{example}
When $L=L_{A, B}$ with $A$ the adding machine and $B$ the Lamplighter, we have
\begin{align*}
\tilde M_0 =
\left(\!
\begin{array}{cc}
 0 & 1 \\
 0 & 0 \\
\end{array}
\!\right)
\hspace{1in}
\tilde M_1 = 
\left(\!
\begin{array}{cc}
 1 & 0 \\
 1 & 1 
\end{array}
\!\right)
\end{align*}
The orbit of $[1:0]$ under the action of $\tilde M_1$ is $\{[1:n]\ :\ n \in \nn\}$, and is not finite.
\tri
\end{example}
\section{When the measure is Markov}

Having obtained the automaton that computes a measure, we can ask the question of what kind of measure it is. Recall from Definition \ref{kstepmarkovmeasuredefinition} that a $k$-step Markov measure is a measure whose sections are uniquely determined by \ii{suffixes} of length $k$, regardless of what comes before. The following theorem provides necessary and sufficient conditions for a finite-state measure to be $k$-step Markov.

\begin{theorem}\label{kstepmark}
Let $\mu$ be a finite-state measure with $n$ nontrivial sections $\mu_1,\ldots,\mu_n$.
Then $\mu$ is $k$-step Markov (for some $k \in \nn$) if and only if for any nonempty word $w \in X^*$, there is at most one $i$, $1\leq i \leq n$ such that $\mu_i|_w = \mu_i$. 
When $\mu$ is $k$-step Markov, we can choose $k\leq n(n-1)$.
\end{theorem}

\pf $\Rightarrow$  Proof by contradiction. Let $\mu$ be $k$-step Markov. Assume that the hypothesis of the theorem does not hold, that is, $\mu$ has two distinct nontrivial sections $\mu|_u$ and $\mu|_v$ (where $u,v\in X^*$) such that $(\mu|_u)|_w = \mu|_u$ and $(\mu|_v)|_w = \mu|_v$ for some nonempty word $w \in X^*$.  Let $W=www\ldots w$ be the word $w$ repeated several times so that $|W|>k$. Then
\begin{align*}
\mu|_{uW} &= (\mu|_u)|_W  = (\mu|_u)|_{ww\ldots w} = \mu|_u,\\
\mu|_{vW} &= (\mu|_v)|_W  = (\mu|_v)|_{ww\ldots w} = \mu|_v.
\end{align*}
So the nontrivial sections $\mu|_{uW}$ and $\mu|_{vW}$ are different, but $|W|>k$. That is, a suffix of length $k$  does not uniquely determine a nontrivial section of $\mu$.  This contradicts the assumption that $\mu$ is $k$-Markov.

$\Leftarrow:$ Assume now that the hypothesis holds. For $\mu$ to be $k$-step Markov, it suffices to show that for any two nontrivial sections $\mu|_u$ and $\mu|_v$, and any word $w$ with $|w|=k$, we have $\mu|_{uw} = \mu|_{vw}$ whenever both $uw$ and $vw$ are admissible words.

Let $k=n(n-1)$ and fix $w=w_1w_2\ldots w_k$, $w_i \in X$ with $|w|=k$. Consider Table \ref{pathstable}. The columns of this table are the paths from $\mu|_u$ and $\mu|_v$ obtained by taking sections by the word $w$ character by character.
\begin{table}[!ht]
\centering
\begin{tabular}{|c|c|}
\hline
$\mu|_u$ & $\mu|_v$\\
$\mu|_{uw_1}$ & $\mu|_{vw_1}$\\
$\mu|_{uw_1w_2}$ & $\mu|_{vw_1w_2}$\\
$\ldots$ & $\ldots$\\
$\mu|_{uw}$ & $\mu|_{vw}$\\
\hline
\end{tabular}
\caption{Paths of length $k$ starting from $\mu|_u$ and $\mu|_v$}
\label{pathstable}
\end{table}
By assumption, $uw$ and $vw$ are admissible, so Table \ref{pathstable} contains only nontrivial sections.  We claim that one row of the table contains two identical measures: $\mu|_{uw_1..w_i}= \mu|_{vw_1..w_i}$ for some $i$, $0\le i\le k$.  Then each subsequent row also contains two identical measures.  In particular, $\mu|_{uw} = \mu|_{vw}$.  To prove the claim (and complete the proof), assume the contrary.  Since $\mu$ has only $n$ nontrivial sections, there are only $n(n-1)$ pairs of distinct nontrivial sections.  Table \ref{pathstable} has $k+1>n(n-1)$ rows, hence a row in the table must repeat, i.e.,
\[
(\mu|_{uw_1..w_i}, \mu|_{vw_1..w_i}) = (\mu|_{uw_1..w_j}, \mu|_{vw_1..w_j})
\]

for some $1 \leq i < j \leq k$. But that means that the word $W:=w_{i+1}\ldots w_j$ fixes two sections $\mu_U := \mu|_{uw_0..w_i}$ and $\mu_V := \mu|_{vw_0..w_i}$; that is, $(\mu|_U)|_W = \mu|_U$ and $(\mu|_V)_W=\mu|_V$. By our hypothesis, the nontrivial section fixed by $W$ is unique, so $\mu|_U = \mu_V$.
Since $UW = uw$ and $VW = vw$, this implies $\mu|_{uw}  = \mu|_{vw}$. This completes the proof. \sqr
\medbreak

\bb{Remark:} The free semigroup $FS(X)$ generated by $X$ acts on the sections of $\mu$: for $w \in FS(X)$, $w \cdot \mu_i := \mu_i|_w$. The condition of Theorem \ref{kstepmark} can be re-stated as follows: the action of any nonempty word $w  \in X^*$ on the sections of $\mu$ has at most one nontrivial fixed point. 

\medbreak

Theorem \ref{kstepmark} is illustrated by the following example.

\begin{example}
When $A$ is the adding machine and $B$ is \alesh, the measure $\mu_{A,B}$ is defined by the automaton $M$ in Figure \ref{fig:aleshmeas}.  The measure satisfies the hypothesis of Theorem \ref{kstepmark}.  By the theorem, $\mu_{A,B}$ is $k$-step Markov for some $k\le20$.  Direct examination of the automaton $M$ reveals that the admissible words for $\mu_{A,B}$ are all words not containing $000$ or $1101$, and the measure is, in fact, $3$-step Markov (see Table \ref{tab:mumarkalesh}).\tri

\end{example}

\begin{table}[!ht]
\centering
\begin{tabular}{|c|c|}
\hline 
$w$ ends in & $\mu|_w$ \\ \hline
$00$ & $\displaystyle\frac{\mu|_b + \mu|_{c,1}}{2}$ \\ 
$11$ & $\displaystyle\frac{\mu|_b + \mu|_{c,A}}{2}$ \\ 
$01$ & $\displaystyle\frac{\mu|_b + \mu|_{c,A} + 2\mu_a}{4}$ \\ 
$110$ & $\displaystyle\mu_a$ \\ 
$010$ & $\displaystyle\frac{\mu_a + \mu|_b + \mu|_{c,1} }{3}$ \\ 
\hline 
\end{tabular} 

Minimal forbidden words: 

$000$, $1101$.
\caption{$\mu_{A,B}$ as a Markov measure when $B$ is \alesh}
\label{tab:mumarkalesh}
\end{table}

The condition of Theorem \ref{kstepmark} is not satisfied trivially.

\begin{figure}[!ht]
\centering
\includegraphics[width=0.5\textwidth]{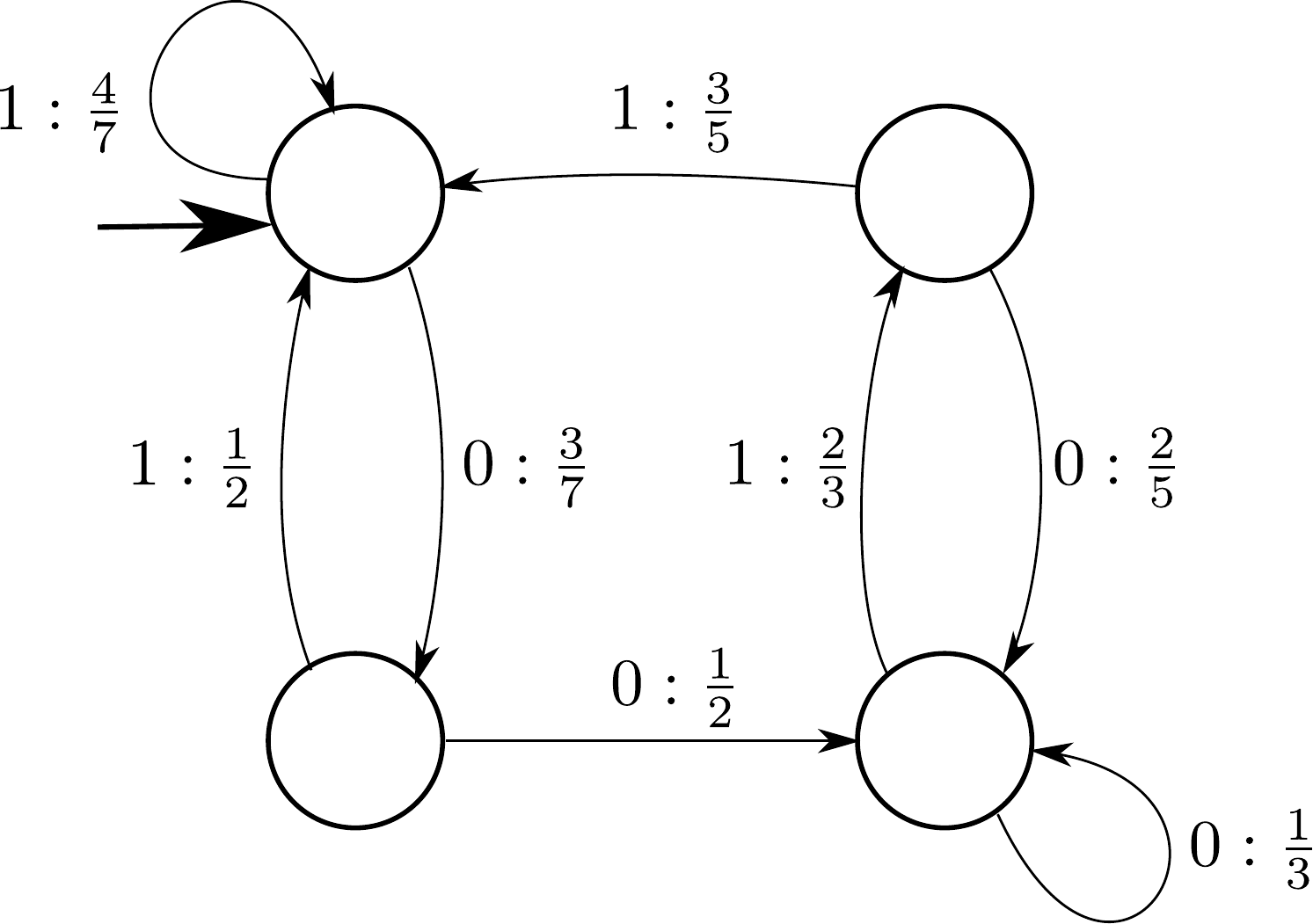}
\caption{A diagram of the automaton of a finite-state measure that is not Markov.}
\label{fig:nonmarkov}
\end{figure}

\begin{example}
A finite-state measure $\mu$ defined by the automaton in Figure \ref{fig:nonmarkov} is not a Markov measure. The initial state $\mu$ (on the top left) is fixed by the action of the word $01$, but so is the top right state, $\mu|_{001}$: $\mu|_{01} = \mu$ and $(\mu|_{001})|_{01} = \mu|_{001}$.
Note that $\mu \neq \mu|_{001}$ since $\mu(0X^\nn)=\frac{3}{7}$ while $\mu|_{001}(0X^\nn) = \frac{2}{5}$. By Theorem \ref{kstepmark}, $\mu$ cannot be a $k$-step Markov measure for any $k$.
\tri
\end{example}

\newpage

\printbibliography

\end{document}